\documentclass[10pt]{article}

\usepackage{a4wide}
\usepackage{amssymb}
\usepackage{amsfonts}
\usepackage{amsmath}
\input xy
\xyoption{arrow} \xyoption{matrix}

\date{}

\newtheorem{proposition}{Proposition}[section]
\newtheorem{theorem}[proposition]{Theorem}
\newtheorem{lemma}[proposition]{Lemma}

\newtheorem{corollary}[proposition]{Corollary}
\newtheorem{question}[proposition]{Question}

\def\der{\partial }

\def\nFM0{{\nu }_{F,M_0}}
\def\nFN0{{\nu }_{F,N_0}}
\def\nGN0{{\nu }_{G,N_0}}

\def\N0{ {\bf N}_0 }

\def\t{\otimes}

\def\v{\varphi}
\def\ra{\rightarrow}

\def\Xpm{X^{\pm }}

\def\s{\sigma}
\def\Z{\mathbb{Z}}

\def\l1{{\lambda}_1}

\def\a{\alpha}
\def\a0{ {\alpha }_0}
\def\a1{ {\alpha }_1}

\def\l{\lambda}
\def\o{\omega}

\def\nFGM0{{\nu }_{F,G,M_0}}

\def\nFN0{{\nu}_{F,N_0}}

\def\sm{{\sigma}^m}

\def\sm1{{\sigma}^{-1}}

\def\smtp1{{\sigma}^{-t+1}}

\def\o{\omega }
\def\S1{S^{-1}}

\def\Xpm1{X^{\pm 1}_1}

\def\sPM1{{\sigma }^{\pm 1}}
\def\sMP1{{\sigma }^{\mp 1 }}

\def\d{\delta}

\def\di{{\rm d.ind}}

\def\L{\Lambda}

\def\Ytm1{Y^{t-1}}
\def\Yim1{Y^{i-1}}

\def\CM{{\cal M}}

\def\CF{{\cal F}}

\def\ass{{\rm ass}}

\def\supp{{\rm supp }}

\def\Aut{{\rm Aut}}

\def\dim{{\rm dim }}

\def\ker{ {\rm ker } }

\def\CJ{ {\cal J}}

\def\SL2Z{ {\rm SL}_2({\bf Z}) }

\def\CR{ {\cal R}}

\def\th{ \theta }

\def\Gp1{ G^{1 , 1 } }
\def\P11{ P^{-1 , 1 } }
\def\Pp1{ P^{1 , 1 } }

\def\Supp{{\rm Supp}}

\def\th{\theta}

\def\nCLsr{{}^\nu\kern-2pt {\cal L}^{\sigma , \rho  }}
\def\nP{{}^\nu \kern-2pt P}
\def\nL{{}^\nu\kern-2pt L}
\def\nLL{{}^\nu\kern-2pt \Lambda}
\def\nPsr{{}^\nu\kern-2pt P^{\sigma , \rho  }}
\def\nLsr{{}^\nu\kern-2pt L^{\sigma , \rho  }}
\def\nuCL{{}^\nu\kern-2pt  {\cal L}}
\def\nCLsr{{}^\nu\kern-2pt {\cal L}^{\sigma , \rho  }}
\def\nCL1m{{}^\nu\kern-2pt {\cal L}^{-1 , 1  }}
\def\x1nu{x^\frac{1}{\nu}}
\def\xm1nu{x^{-\frac{1}{\nu}}}

\def\CR{ {\cal R}}

\def\ra{\rightarrow }

\def\CB{{\cal B}}

\def\coker{{\rm coker}}
\def\CT{{\cal T}}

\def\CC{ {\cal C}}

\def\nAM0{{\nu }_{{\cal A},M_0}}
\def\nAN0{{\nu }_{{\cal A},N_0}}

\def\End{ {\rm End }}

\def\CJ{ {\cal J }}
\def\CR{ {\cal R }}

\def\bR{\overline{R}}

\def\ga{\mathfrak{a}}
\def\gb{\mathfrak{b}}

\def\gm{\mathfrak{m}}
\def\gp{\mathfrak{p}}

\def\SL{{\rm SL}}

\def\di!{\frac{\der^i}{i!}}
\def\dik!{\frac{\der^k_i}{k!}}

\def\N{\mathbb{N}}

\def\0{\overline{0}}
\def\1{\overline{1}}

\def\Ln1{\L_{n,\overline{1}}}

\def\oa{\overline{a}}

\def\a1{a_{\overline{1}}}

\def\bs{\overline{s}}

\def\S{\Sigma}

\def\vn1{\overrightarrow{n-1}}

\def\im{{\rm im}}

\def\Min{{\rm Min}}

\def\Inn{{\rm Inn}}

\def\mJ{\mathbb{J}}
\def\mI{\mathbb{I}}

\def\ind{{\rm ind}}

\def\Frac{{\rm Frac}}

\def\K1{{\rm K}_1}

\def\hmI1{\widehat{\mI_1}}
\def\tmI1{\widetilde{\mI_1}}
\def\tmJ1{\widetilde{\mJ_1}}
\def\hB1{\widehat{B_1}}
\def\hCB1{\widehat{\CB_1}}

\def\bS{\overline{S}}

\def\Den{{\rm Den}}
\def\Denl{{\rm Den}_l}
\def\Ore{{\rm Ore}}
\def\ore{{\rm ore}}
\def\Den{{\rm Den}}
\def\LDen{{\rm LDen}}
\def\minLDen{{\rm min.LDen}}
\def\Loc{{\rm Loc}}
\def\minLoc{{\rm min.Loc}}
\def\Ass{{\rm Ass}}
\def\ttS{\widetilde{S}}
\def\ot{\overline{t}}
\def\maxDen{{\rm max.Den}}
\def\maxAss{{\rm max.Ass}}
\def\maxLoc{{\rm max.Loc}}
\def\llrad{{\rm l.lrad}}
\def\mW{\mathrm{W}}
\def\assmaxDen{{\rm ass.max.Den}}

\begin{document}

\author{V. V. \  Bavula 
}

\title{The largest left quotient ring of a ring}

\maketitle

\begin{abstract}

The left quotient ring (i.e. the left classical ring of fractions)
$Q_{cl}(R)$ of a ring $R$ does not always exist and still, in
general, there is no good understanding of the reason why this
happens. In this paper, it is proved existence of the {\em largest
left quotient ring} $Q_l(R)$, i.e. $Q_l(R) = S_0(R)^{-1}R$ where
$S_0(R)$ is the {\em largest left regular denominator set} of $R$.
It is proved that $Q_l(Q_l(R))=Q_l(R)$; the ring $Q_l(R)$ is
semi-simple iff $Q_{cl}(R)$ exists and is semi-simple; moreover,
if the ring $Q_l(R)$ is left artinian then $Q_{cl}(R)$ exists and
$Q_l(R) = Q_{cl}(R)$. The group of units $Q_l(R)^*$ of $Q_l(R)$ is
equal to the set $\{ s^{-1} t\, | \, s,t\in S_0(R)\}$ and $S_0(R)
= R\cap Q_l(R)^*$. If there exists a  finitely generated flat left
$R$-module which is not projective then $Q_l(R)$ is not a
semi-simple ring. We extend slightly Ore's method of localization
to {\em localizable left Ore sets}, give a criterion of when a
left Ore set is localizable, and prove that {\em all left and
right} Ore sets of an arbitrary ring  are localizable (not just
denominator sets as in Ore's method of localization). Applications
are given for certain classes of rings (semi-prime Goldie rings,
 Noetherian commutative  rings, the algebra of polynomial integro-differential operators).


{\em Key Words: the largest left quotient ring of a ring, the
largest left (regular) denominator set of a ring,
  the classical left quotient ring of a ring,   denominator set,
   the maximal left quotient rings of a
ring. }

 {\em Mathematics subject classification
2000: 16U20,  16P40, 16S32, 13N10.}

$${\bf Contents}$$
\begin{enumerate}
\item Introduction. \item The largest left quotient ring of a
ring. \item The maximal left quotient rings of a ring. \item The
largest  quotient ring of a ring.
\item Examples.
\end{enumerate}
\end{abstract}


\section{Introduction}

The aim of the paper is to introduce the following new concepts
and to prove their existence for an {\em arbitrary} ring: {\em the
largest left quotient ring of a ring, the largest left (regular)
denominator set of a ring,  the maximal left quotient rings of a
ring, the largest (two-sided) quotient ring of a ring, the maximal
(two-sided) quotient rings of a ring, a (left) localization
maximal ring}.

Throughout, module means a left module; $R$ is a ring with 1;
$\CC_R$ is the set of (left and right) regular elements of the
ring $R$  (i.e. $\CC_R$ is the set of non-zero-divisors of $R$);
$Q_{cl}(R)$ is the left quotient ring of $R$ (if exists); $\Denl
(R,0)$ is the set of regular  left denominator sets $S$ in $R$
($S\subseteq \CC_R$).

$\noindent $

{\bf The largest left regular denominator set and the largest left
quotient ring of a ring}.

\begin{itemize}
\item (Theorem \ref{3Jul10}) {\em There exists the largest (w.r.t.
 inclusion) left regular   denominators set $S_0(R)$ in $R$ and so
 $Q_l(R):= S_0(R)^{-1}R$ is the largest left quotient ring of}  $R$.
 \item (Corollary \ref{a7Jul10}) {\em The set $\Denl (R,0)$ of
 left regular   denominator sets of $R$ is a complete lattice and an
 abelian monoid.}
\end{itemize}

The next theorem describes the group  of units $Q_l(R)^*$ of the
ring $Q_l(R)$  and its connection with $S_0(R)$.

\begin{itemize} \item (Theorem \ref{4Jul10})

\begin{enumerate}
\item $ S_0 (Q_l(R))= Q_l(R)^*$ {\em and} $S_0(Q_l(R))\cap R=
S_0(R)$.
 \item $Q_l(R)^*= \langle S_0(R), S_0(R)^{-1}\rangle$, {\em i.e. the
 group of units of the ring $Q_l(R)$ is generated by the sets
 $S_0(R)$ and} $S_0(R)^{-1}:= \{ s^{-1} \, | \, s\in S_0(R)\}$.
 \item $Q_l(R)^* = \{ s^{-1}t\, | \, s,t\in S_0(R)\}$.
 \item $Q_l(Q_l(R))=Q_l(R)$.
\end{enumerate}
\end{itemize}

The next theorem gives an answer to the question of when the ring
$Q_l(R)$ is semi-simple (Theorem \ref{4Jul10} is a key step in
proving Theorem \ref{5Jul10}, statements 2-5 is Goldie's Theorem).

\begin{itemize} \item (Theorem \ref{5Jul10}) {\em The following
statements are equivalent.}

\begin{enumerate}
\item $Q_l(R)$ {\em is a semi-simple ring.} \item $Q_{cl}(R)$ {\em
exists and is a semi-simple ring.}
 \item $R$ {\em  is a left order in a semi-simple ring.}
 \item $R$ {\em has finite left rank, satisfies the ascending chain
 condition on left annihilators, and is a semi-prime ring.}
 \item {\em A left ideal of the ring $R$ is essential iff it contains a
 regular element.}
\end{enumerate}
{\em If one of these conditions hold then $S_0(R) = \CC_R$ and}
$Q_l(R) = Q_{cl}(R)$.
\end{itemize}

The next result is a sufficient condition for the ring $Q_l(R)$
 not being a semi-simple ring.

\begin{itemize}
\item (Corollary \ref{11Jul10}) {\em If there exists a finitely
generated flat $R$-module  which is not projective then the ring
$Q_l(R)$ is not semi-simple.}
\end{itemize}

The next corollary gives a sufficient condition for existing of
$Q_{cl}(R)$.

\begin{itemize}
\item (Corollary \ref{A5Jul10}) {\em  If $Q_l(R)$ is a left
artinian ring then $S_0(R) = \CC_R$ and $Q_l(R) = Q_{cl}(R)$.}
\end{itemize}

Let $S_0^r(R)$ and $Q_r(R):=RS_0^r(R)^{-1}$ be the {\em  largest
regular right denominator set} in $R$ and the {\em  largest right
quotient ring} of $R$, respectively. In general, the following
natural questions have negative answers as the algebra $\mI_1:=
K\langle x, \frac{d}{dx}, \int\rangle$ of polynomial
integro-differential operators over a field $K$ of characteristic
zero demonstrates \cite{intdifline}:

{\it Question 1}. {\em Is $S_0(R)=S_0^r(R)$?}

{\it Question 2}. {\em Is $S_0(R)\subseteq S_0^r(R)$ or
$S_0(R)\supseteq S_0^r(R)$?}

{\it Question 3}. {\em  Are the rings $Q_l(R)$ and $Q_r(R)$
isomorphic?}

Though, for the algebra $\mI_1$ the next question has positive
answer \cite{intdifline}.

{\it Question 4}. {\em  Are the rings $Q_l(R)$ and $Q_r(R)$
anti-isomorphic?}

{\it Remark}. The algebra $\mI_1$ is neither left nor right
Noetherian, it contains infinite direct sums of nonzero left (and
right) ideals, neither left nor right quotient ring of $\mI_1$
exists \cite{intdifline}.

$\noindent $

{\bf Notation}:

\begin{itemize}
\item $\Ore_l(R):=\{ S\, | \, S$ is a left Ore set in $R\}$; \item
$\Den_l(R):=\{ S\, | \, S$ is a left denominator set in $R\}$;
\item $\Loc_l(R):= \{ S^{-1}R\, | \, S\in \Den_l(R)\}$; \item
$\Ass_l(R):= \{ \ass (S)\, | \, S\in \Den_l(R)\}$ where $\ass
(S):= \{ r\in R \, | \, sr=0$ for some $s=s(r)\in S\}$; \item
$\Den_l(R, \ga ) := \{ S\in \Den_l(R)\, | \, \ass (S)=\ga \}$
where $\ga \in \Ass_l(R)$; \item $S_\ga=S_\ga (R)=S_{l,\ga }(R)$
is the {\em largest element} of the poset $(\Den_l(R, \ga ),
\subseteq )$ and $Q_\ga (R):=Q_{l,\ga }(R):=S_\ga^{-1} R$ is  the
{\em largest left quotient ring associated to} $\ga$, $S_\ga $
exists (Theorem \ref{3Jul10}.(2)); \item In particular,
$S_0=S_0(R)=S_{l,0}(R)$ is the largest element of the poset
$(\Den_l(R, 0), \subseteq )$ and $Q_l(R):=S_0^{-1}R$ is the
largest left quotient ring of $R$; \item $\Loc_l(R, \ga ):= \{
S^{-1}R\, | \, S\in \Den_l(R, \ga )\}$.
\end{itemize}

$\noindent $

For each denominator set $S\in \Den_l(R, \ga )$ where $\ga := \ass
(S)$, there are natural ring homomorphisms $R\stackrel{\pi}{\ra}
R/ \ga \ra S^{-1}R$. In Section \ref{TMLQR}, connections between
the denominator sets $\Den_l(R, \ga )$, $\Den_l( R/ \ga , 0)$ and
$\Den_l(S^{-1}R, 0)$ are established (Lemma \ref{a5Jul10}, Lemma
\ref{a6Jul10} and Proposition \ref{8Jul10}; these results  are too
technical to explain in the Introduction). They are used to prove
the following results.
\begin{itemize}
\item (Lemma \ref{a6Jul10}.(2)) {\em Let $a\in \Ass_l(R)$ and $\pi
: R\ra R/ \ga$, $a\mapsto a+\ga$. Then $\pi^{-1} (S_0(R/ \ga )) =
S_\ga (R)$, $\pi (S_\ga (R)) = S_0(R/ \ga )$ and } $ Q_\ga (R) =
Q_l( R/ \ga )$. \item (Lemma \ref{b14Nov10}.(2)) {\em The set
$\maxDen_l(R)$ of maximal elements of the poset $(\Den_l(R),
\subseteq )$  of left denominator sets of an arbitrary  ring $R$
is a non-empty set.}
\end{itemize}
This means that the set $\maxLoc_l(R)$ of maximal left quotient
rings of an arbitrary ring $R$ is a non-empty set, and
$\maxLoc_l(R)=\{ S^{-1}R\, | \, S\in \maxDen_l(R)\}$.

$\noindent $

{\it Example}. $\maxDen_l(\mI_1 ) = \{ \mI_1\backslash F\}$
(Proposition \ref{a11Dec10}.(3)) where $F$ is the only proper
ideal of the algebra of polynomial integro-differential operators
$\mI_1$ and $\maxLoc_l(\mI_1) = \{ \Frac (A_1)\}$ (Proposition
\ref{a11Dec10}.(2)) where $\Frac (A_1)$ is the {\em Weyl skew
field}, i.e. the classical (left and right)  ring of fractions of
the first Weyl algebra $A_1= K \langle x,
\frac{x}{dx}\rangle$.  

$\noindent $

A new class of rings is introduced, the class of {\em left
localization maximal rings}, see Section \ref{TMLQR} (intuitively,
one can invert nothing more in these rings). For an arbitrary ring
$R$, a criterion is given in terms of this class of rings of when
a left localization $S^{-1}R$ of $R$ is a maximal left
localization of $R$.
\begin{itemize}
\item (Theorem \ref{21Nov10}) {\em Let  a ring $A$ be a left
localization of a ring $R$, i.e. $A\in \Loc_l(R, \ga )$ for some
$\ga \in \Ass_l( R)$. Then $A\in \maxLoc_l( R)$ iff $Q_l( A) = A$
and  $\Ass_l(A) = \{ 0\}$, i.e. $A$ is a left localization maximal
ring. }
\end{itemize}

{\it Example}. Let $V$ be an infinite dimensional vector space
with  countable basis over a field $K$, let $\CC :=\{ \v \in
\End_K(V)\, | \, \dim_K(\im (\v ))<\infty \}$ be the ideal of
compact operators/linear maps of the algebra $\End_K(V)$. Then the
factor algebra $\End_K(V) / \CC$ is a left localization maximal
ring (Theorem \ref{4Dec10}.(5)). Moreover, the set $\CF$ of
Fredholm linear maps in $V$ is the maximal left  (resp. right;
left and right) denominator set of the ring $\End_K(V)$ and $$
\CF^{-1}\End_K(V)\simeq \End_K(V)\CF^{-1} \simeq \End_K(V) / \CC$$
is the maximal left (resp. right; left and right) localization
ring of $\End_K(V)$ (Theorem \ref{4Dec10}).

{\bf A slight extension/generalization of Ore's method of
localization}. Ore's method of left localization says that we can
localize a ring $R$ {\em precisely} at the left denominator sets
$\Den_l(R)$ of the ring $R$. Notice that each left denominator set
is a left Ore set but not the other way round, in general. We
introduce the concept of a {\em localizable left Ore set} and give
a criterion of when a left Ore set is localizable (Theorem
\ref{A15Nov10}). Each left denominator set is a localizable left
Ore set but not vice versa, in general. We extend Ore's method of
localization to localizable left Ore sets and prove an analogue of
Ore's Theorem (by using Ore's Theorem) for localizable left Ore
sets.

\begin{itemize}
\item (Corollary \ref{28Nov10}) {\em Let $S$ be a localizable left
Ore set in a ring $R$. Then there exists an ordered pair $(Q, f)$
where $Q$ is a ring and $f: R\ra Q$ is a ring homomorphism such
that

(i)  for all $s\in S$, $f(s)$ is a unit in $Q$; 


 and if $(Q', f')$ is another
 pair  satisfying the condition
 (i) then there is a unique ring homomorphism $h: Q\ra Q'$ such that
$f' = hf$. The ring $Q$ is unique up to isomorphism. The ring $Q$
is isomorphic to the left localization of the ring $R/ \gp (S)$ at
the left denominator set $\pi (S)\in \Den_l(R/ \gp (S), 0)$ where
the ideal $\gp (S)$ of the ring $R$ is defined in (\ref{paup2})
and $\pi : R\ra R/ \gp (S)$, $ a\mapsto a+\gp (S)$. }
\end{itemize}

In Section \ref{LQR}, we prove {\em two-sided} analogues of some
of the results of Sections \ref{LLQRR} and \ref{TMLQR}. In most
cases the proofs are easy corollaries of the corresponding
one-sided (i.e. left and right) results but there are some
surprises. In particular,
\begin{itemize}
\item (Theorem \ref{B15Nov10}) {\em Every (left and right) Ore set
is localizable.}
\end{itemize}
Therefore we can localize at {\em all} the (left and right) Ore
sets not just at the  (left and right) denominator sets as in
Ore's method of localization.

\begin{itemize}
\item (Corollary \ref{T28Nov10}) {\em Let $S$ be an  Ore set in a
ring $R$. Then there exists an ordered pair $(Q, f)$ where $Q$ is
a ring and $f: R\ra Q$ is a ring homomorphism such that

(i)  for all $s\in S$, $f(s)$ is a unit in $Q$; 


 and if $(Q', f')$ is another
 pair satisfying the condition
 (i) then there is a unique ring homomorphism $h: Q\ra Q'$ such that
$f' = hf$. The ring $Q$ is unique up to isomorphism. The ring $Q$
is isomorphic to the  localization of the ring $R/ \gp (S)$ at the
 denominator set $\pi (S)\in \Den (R/ \gp (S), 0)$ where the ideal
$\gp (S)$ of the ring $R$ is defined in (\ref{paup4}) and $\pi :
R\ra R/ \gp (S)$, $ a\mapsto a+\gp (S)$. }
\end{itemize}

 In Section \ref{EXEM}, the largest (left; right; left and right)
and maximal (left; right; left and right) quotient rings are found
for following rings: the endomorphism algebra $\End_K(V)$ of an
infinite dimensional vector space $V$ with countable basis,
semi-prime Goldie rings,  the  algebra $\mI_1$ of polynomial
integro-differential operators and  Noetherian commutative rings.


\section{The largest left quotient ring of a ring}\label{LLQRR}

In this section, existence of the largest left quotient ring of a
ring is proved (Theorem \ref{3Jul10}). Proofs of Theorems
\ref{4Jul10} and \ref{5Jul10} are given.

$\noindent $

{\bf The largest left quotient ring of a ring}. Let $R$ be a ring.
A {\em multiplicatively closed subset} $S$ of $R$ (i.e. a
multiplicative sub-semigroup of $(R, \cdot )$ such that $1\in S$
and $0\not\in S$) is said to be a {\em left Ore set} if it
satisfies the {\em left Ore condition}: for each $r\in R$ and
$s\in S$, $Sr\bigcap Rs\neq \emptyset$. Let $S$ be a (non-empty)
multiplicatively closed subset of $R$, and let $\ass (S) :=\{ r\in
R\, | \, sr=0$ for some $s\in S\}$ (if, in addition, $S$ is a left
Ore set then $\ass (S)$ is an ideal of the ring $R$). Then a {\em
left quotient ring} of $R$ with respect to $S$ (a {\em left
localization} of $R$ at $S$) is a ring $Q$ together with a
homomorphism $\v :R\ra Q$ such that

(i) for all $s\in S$, $\v (s)$ is a unit of $Q$,

(ii) for all $q\in Q$, $q=\v (s)^{-1}\v (r) $ for some $r\in R$,
$s\in S$, and

(iii) $\ker (\v ) = \ass (S)$.

If exists the ring $Q$ is unique up to isomorphism, usually it is
denoted by $S^{-1}R$. Recall that $S^{-1}R$ exists iff $S$ is a
left Ore set and the set $\bS=\{ s+\ass (S)\in R/\ass (S)\, | \,
s\in S\} $  consists of regular elements (\cite{MR}, 2.1.12). If
the last two conditions are satisfied then $S$ is called a {\em
left denominator set}.  Similarly, a {\em right Ore set}, the {\em
right Ore condition}, the {\em right denominator set}  and the
{\em right quotient ring} $RS^{-1}$ are defined. If both $S^{-1}R$
and $RS^{-1}$ exist then they are isomorphic (\cite{MR},
2.1.4.(ii)). The left quotient ring of $R$ with respect to the set
$\CC_R$ of all regular elements is called the {\em left quotient
ring} of $R$, if exists, it is denoted by $\Frac_l(R)$ or
$Q_{cl}(R)$. Similarly, the {\em right quotient ring},
$\Frac_r(R)=Q_{cl}^r(R)$, is defined. If both left and right
quotient rings of $R$ exist then they are isomorphic and we write
simply $\Frac (R)$ or $Q(R)$ in this case.


\begin{theorem}\label{3Jul10}
\begin{enumerate}
\item For each $\ga \in \Ass_l(R)$, the set $\Den_l(R, \ga )$   is
an ordered abelian semigroup $(S_1S_2= S_2S_1$, and $S_1\subseteq
S_2$ implies $S_1S_3\subseteq S_2S_3$) where the product $S_1S_2=
\langle S_1, S_2\rangle$ is the multiplicative subsemigroup of
$(R, \cdot )$ generated by $S_1$ and $S_2$.  \item  $S_\ga :=
S_\ga (R) := \bigcup_{S\in \Den_l(R, \ga )}S$ is the largest
element (w.r.t. $\subseteq $) in $\Den_l(R, \ga )$. The set
$S_\ga$ is called the {\em largest left denominator set}
associated to $\ga$.   \item Let $S_i\in \Den_l(R,\ga )$, $i\in
I$, where $I$ is an arbitrary non-empty set. Then
\begin{equation}\label{SiIgen}
\langle S_i\, | \, i\in I\rangle:=\bigcup_{\emptyset \neq
J\subseteq I, |J|<\infty}\prod_{j\in J}S_j\in \Den_l(R, \ga )
\end{equation}
the left denominator set generated by the left denominators sets
$S_i$, it is the least upper bound of the set $\{ S_i\}_{i\in I}$
in $\Den_l(R, \ga )$, i.e. $\langle S_i\, | \, i\in I\rangle=
\bigvee_{i\in I}S_i$.
\end{enumerate}
\end{theorem}

{\it Remark}. Clearly, $S_1S_2= S_1\bigvee S_2$, the {\em join} of
$S_1$ and $S_2$.

{\it Proof}. 1. It remains to show that $S_1S_2\in \Den_l(R, \ga
)$, that is, $S_1S_2$ is a left Ore set in $R$ with $\ass (S_1S_2)
= \ga $ and $rs=0$ where $r\in R$ and $s\in S$ implies $tr=0$ for
some $t\in S$. To prove that the left Ore condition holds for the
multiplicatively closed set $S=S_1S_2$ we have to show that, for
each $s\in S$ and $ a\in R$, there exist elements $t\in S$ and
$b\in R$ such that $ta = bs$. We use induction on the length $l(s)
= \min \{ n \, | \, s=s_1\cdots s_n$ where all $s_i\in S_1\cup
S_2\}$ of the element $s$. If $l(s) =1$, i.e. $s\in S_1\cup S_2$,
then the statement is obvious since $S_1$ and $S_2$ are left Ore
sets. Suppose that $n:= l(s)>1$, and the statement holds for
$n'<n$. Fix a presentation $s=s_1\cdots s_{n-1}s_n$ where all
$s_i\in S_1\cup S_2$. Then $t_1a=b_1s_n$ for some elements $t_1\in
S_1\cup S_2$ and $b_1\in R$, and, by induction, $t_2b_1=
bs_1\cdots s_{n-1}$ for some elements $t_2\in S$ and $b\in R$.
Clearly, $t=t_2t_1 \in S$ and $ta=t_2(t_1a)= (t_2b_1) s_n= bs$, as
required.

Let us show that  $\ass (S)=\ga$.  Clearly, $\ass (S_1S_2)
\supseteq \ga$. Let $r\in \ass (S_1S_2)$. Then $sr=0$ for some
element $s=s_1\cdots s_m\in S$ where all $s_i\in S_1\cup S_2$.
Without loss of generality we may assume that $s_{odd}\in S_1$ and
$s_{even}\in S_2$. Then $s_2\cdots s_m r\in \ga$, and so there
exists $s_2'\in S_2$ such that $(s_2's_2)s_3\cdots s_mr=
s_2''s_3\cdots s_mr=0$ where $s_2''= s_2's_2\in S_2$. By induction
on $m$, we conclude that $r\in \ga$. Therefore, $\ass (S) = \ga$.

Finally, suppose that $rs=0$ where $r\in R$  and $s\in S$. We have
to show that $tr=0$ for some $t\in S$. Fix a presentation
$s=s_1\ldots s_n$ where $s_i\in S_1\cup S_2$. Then
$0=rs=rs_1\cdots s_n$ implies that $t_1rs_1\cdots s_{n-1}=0$ for
some $t_1\in S_1\cup S_2$ since $S_1, S_2\in \Den_l(R,\ga )$.
Similarly, $t_2t_1rs_1\cdots s_{n-2}=0$ for some $t_2\in S_1\cup
S_2$. Repeating the same argument several times we see that $t_n
t_{n-1}\cdots t_1r=0$ for some $t_i\in S_1\cup S_2$. Notice that
$t:=t_n t_{n-1}\cdots t_1\in S$ and $tr=0$, as required.

2. Let $s_1, s_2\in S_\ga$. Then $s_1\in S_1$ and $s_2\in S_2$ for
some $S_1, S_2\in \Den_l(R, \ga )$, and so $s_1,s_2\in S_1S_2\in
\Den_l(R, \ga )$, by statement 1. Therefore, $S_\ga\in
\Den_l(R,\ga )$, and so $S_\ga$ is the largest element in
$\Den_l(R, \ga )$.

3. Statement 3 follows from statement 1.  $\Box $

$\noindent $

{\it Definition}. For each ideal $\ga \in \Ass_l(R)$, the ring
$Q_\ga (R):= S_\ga^{-1}R$ is called the {\em largest left quotient
ring associated to} $\ga$. When $\ga =0$, the ring $Q_l(R):= Q_\ga
(R) := S_0^{-1}R$ is called the {\em largest left quotient ring}
of $R$ and $S_0=S_0(R)$ is called the {\em largest left regular
denominator set} of $R$.

$\noindent $

The next obvious corollary shows that $Q_l(R)$ is a generalization
of the classical left quotient ring $Q_{cl}(R)$.

\begin{corollary}\label{a3Jul10}
\begin{enumerate}
\item If the classical left quotient ring $Q_{cl}(R):=\CC_R^{-1}R$
exists then the set of regular elements $\CC_R$  of the ring $R$
is the largest left regular   denominator set  and
$Q_l(R)=Q_{cl}(R)$.  \item Let $R_1, \ldots ,R_n$ be rings. Then
$Q_l(\prod_{i=1}^n R_i)\simeq \prod_{i=1}^n Q_l(R_i)$.
\end{enumerate}
\end{corollary}

{\it Proof}. It is obvious.  $\Box $

\begin{question}\label{fAB}
Can a ring monomorphism $f: A\ra B$ be lifted (necessarily
uniquely) to a ring monomorphism $f':Q_l(A)\ra Q_l(B)$?
\end{question}

In general, the answer is no.

\begin{proposition}\label{21Aug10}
{\rm (Corollary 9.9, \cite{intdifline})} Let $K$ be a field of
characteristic zero and $A_1=K\langle x, \frac{d}{dx}\rangle$ be
the ring of polynomial differential operators (the first Weyl
algebra). Then the inclusion $A_1\ra \mI_1$ cannot be lifted
neither to a ring homomorphism $\Frac (A_1)=Q_l(A_1)\ra
Q_l(\mI_1)$ nor to $\Frac (A_1)=Q_r(A_1)\ra Q_r(\mI_1)$.
\end{proposition}

Let $(P, \leq )$ be a {\em partially ordered set}, a {\em poset},
for short ($a\leq a$ for all $a\in P$; $a\leq b$ and $b\leq a$
implies $a=b$;  $a\leq b$ and $b\leq c$ implies $a\leq c$). For a
subset $S$ of $P$, an element $x\in P$ is called an {\em upper
bound} (resp. a {\em lower bound}) if $s\leq x$ (resp. $ s\geq x$)
for all $s\in S$. The {\em least upper bound} for $S$ is the least
element in the the set of all upper bounds for $S$. Similarly, the
{\em greatest lower bound}  for $S$ is defined. A {\em lattice} is
a  poset such that each pair $x$, $y$ of elements of the set $S$
has both the least upper bound $x\bigvee y$ (called the {\em join}
for $x$ and $y$) and the greatest lower bound $x\bigwedge y$
(called the {\em meet} of $x$ and $y$). It follows then by
induction that every non-empty finite set of elements has the join
and the meet. A lattice $L$ is {\em complete} if every subset $S$
of $L$ has the {\em least upper bound} (the {\em join} of $S$)
denoted $\sup (S)$ or $\bigvee_{s\in S} s$ and the {\em greatest
lower bound} (the {\em meet} of $S$) denoted $\inf (S)$ or
$\bigwedge_{s\in S}s$. In a complete lattice there exists the
greatest element $\sup (L)$ denoted by 1, and the smallest element
$\inf (L)$ denoted by 0. By definition, $\sup (\emptyset )= 0$.
Let $(L, \leq )$ be a lattice and $a,b\in L$ with $a\leq b$. The
set $[a,b]:= \{ x\in L\, | \, a\leq x\leq b\}$ is the {\em
interval} between $a$ and $b$. The interval $[a,b]$ is a
sublattice of $L$ with $\inf ([a,b])=a$ and $ \sup ([a,b])= b$.


\begin{corollary}\label{a7Jul10}
The abelian monoid $\Den_l(R, 0)$ is a complete lattice such that
 $S_1S_2=S_1\bigvee S_1$ and $\bigwedge_{i\in I}S_i=\bigvee_{S\in
 \Den_l(R, 0), S\subseteq \cap_{i\in I}S_i}S_i$ where all $S_i\in
 \Den_l(R, 0)$.
\end{corollary}

{\it Proof}. The largest and least  elements of the poset
$\Den_l(R, 0)$ are $S_0(R)$ and $\{ 1\}$ respectively. Clearly,
$\{ 1\}$ is the identity element of the semigroup $\Den_l(R, 0)$.
By Theorem \ref{3Jul10}.(3), each nonempty subset of $\Den_l(R,
0)$ has the least upper bound, hence $\Den_l(R, 0)$ is a complete
lattice by Proposition 1.2, Sect. 3, \cite{Stenstrom-RingsQuot}
and $\bigwedge_{i\in I}S_i=\bigvee_{S\in
 \Den_l(R, 0), S\subseteq \cap_{i\in I}S_i}S_i$. $\Box $

$\noindent $

Clearly, $\bigwedge_{i\in I}S_i$ is the largest element of the set
$\{ S\, | \, S\in \Den_l(R,0), S\subseteq \bigcap_{i\in I}S_i\}$.

\begin{corollary}\label{7Jul10}
\begin{enumerate}
\item Let $R$ be a ring. Each ring automorphism $\s\in \Aut(R)$ of
the ring $R$ has the unique extension $\s\in \Aut(Q_l(R))$ to an
automorphism of the ring $Q_l(R)$ given by the rule $\s (s^{-1} r)
= \s (s)^{-1}\s (r)$ where $s\in S_0(R)$ and $r\in R$. \item The
group $\Aut(R)$ is a subgroup of the group $\Aut(Q_l(R))$.
Moreover, $\Aut (R) = \{ \tau \in \Aut(Q_l(R))\, | \, \tau (R) =
R\}$.
\end{enumerate}
\end{corollary}

{\it Proof}. 1. By the uniqueness of the set $S_0(R)$, $\s
(S_0(R)) = S_0(R)$ for all elements $\s \in \Aut (R)$. Now,
statement 1 follows from the universal property of the
localization at $S_0(R)$.

2. Statement 2 follows from statement 1. $\Box $

$\noindent $

For a ring $R$, let $R^*$ be its group of units and $\Inn (R):=\{
\o_u\, | \, u\in R^*\}$ be the {\em group of inner automorphisms}
of the ring $R$ where $\o_u (r) = uru^{-1}$ is the {\em inner
automorphism} determined by the element $u$. The next proposition
is used in the proof of Theorem \ref{4Jul10}.

\begin{proposition}\label{A4Jul10}
Let $R$ be a ring, $S\in \Den_l(R, 0)$, and $T\in \Den_l(S^{-1}R,
0)$; and so $R\subseteq S^{-1}R\subseteq T^{-1} (S^{-1}R)$ are
natural inclusions of rings. Then
\begin{enumerate}
\item $ T^{-1} (S^{-1}R)= T_1^{-1} (S^{-1}R)$ for some $T_1\in
\Den_l(S^{-1}R, 0)$ such that $S, S^{-1}\subseteq T_1$; in
particular, $sT_1s^{-1} = T_1$ for all $s\in S$. \item If, in
addition, $S\subseteq T$ then $T\cap R\in \Den_l(R, 0)$ and
$S\subseteq T\cap R$.
\end{enumerate}
\end{proposition}

{\it Proof}. 1. For each $s\in S$, $T^{-1}(S^{-1}R) = \o_s
(T^{-1}(S^{-1}R))= \o_s(T)^{-1} (S^{-1}R)$ and $\o_s (T) \in
\Den_l(S^{-1}R, 0)$. It suffices to take $T_1:= \langle S, S^{-1},
\o_s (T)\, | \, s\in S\rangle$ in $\Den_l(S^{-1}R, 0)$ since
clearly $S,S^{-1}\in \Den_l(S^{-1}R,0)$ and, for each non-empty
finite subset $J$ of $S$, $(\prod_{s\in J}\o_s(T))^{-1} (S^{-1}R)
= T^{-1}(S^{-1}R)$. In more detail,
$$ T_1^{-1}(S^{-1}R) = \bigcup_{\emptyset \neq J\subseteq I,
|J|<\infty}(\prod_{s\in J}\o_s(T))^{-1} (S^{-1}R)
=\bigcup_{\emptyset \neq J\subseteq I, |J|<\infty}
T^{-1}(S^{-1}R)= T^{-1}(S^{-1}R).$$

2. The set $T':= T\cap R$ is a multiplicatively closed subset of
$\CC_R$ that contains the set $S$. It remains to show that the
left Ore condition holds for $T'$ in $R$: for each elements $t'\in
T'$ and $r\in R$, $T'r\cap Rt'\neq \emptyset$. Since $T\in
\Den_l(S^{-1}R, 0)$, $Tr\cap S^{-1}Rt'\neq \emptyset$. Take an
element, say $u$, from the intersection. The element $u$ can be
written in two different ways as follows
$$ s_1^{-1} t_1\cdot r = s_2^{-1}r_2t'$$
for some $s_1^{-1}t_1\in T$ and $ s_2^{-1}r_2\in S^{-1}R$ where
$s_1, s_2\in S$ and $t_1, r_2\in R$. Clearly, $t_1=s_1\cdot
s_1^{-1} t_1\in R\cap T = T'$  (since $S\subseteq T$), and
$s_1s_2^{-1} = s_3^{-1} r_3$ for some $s_3\in S$ and $r_3\in R$.
Then the element
$$ s_3t_1\cdot r= s_3s_1\cdot s_1^{-1} t_1r=s_3s_1s_2^{-1}r_2t'=
r_3r_2t'$$ belongs to the intersection $ Tr\cap Rt'$ since
$s_3t_1\in T$ (since $S\subseteq T$) and $r_3r_2\in R$.  $\Box $

$\noindent $

{\bf The group of units $Q_l(R)^*$ of $Q_l(R)$}. For a ring $R$
and its largest left quotient ring $Q_l(R)$, Theorem \ref{4Jul10}
is used in the proof of Theorem \ref{5Jul10} and  gives an answers
to the following natural questions:
\begin{itemize}
\item  {\em What is $S_0(Q_l(R))$?} \item  {\em What is
$S_0(Q_l(R))\cap R$?} \item {\em What is the group $Q_l(R)^*$ of
units of the ring $Q_l(R)$?}\item {\em Is the natural inclusion
$Q_l(R)\subseteq Q_l(Q_l(R))$ an equality?}
\end{itemize}

\begin{theorem}\label{4Jul10}
\begin{enumerate}
\item $ S_0 (Q_l(R))= Q_l(R)^*$ {\em and} $S_0(Q_l(R))\cap R=
S_0(R)$.
 \item $Q_l(R)^*= \langle S_0(R), S_0(R)^{-1}\rangle$, {\em i.e. the
 group of units of the ring $Q_l(R)$ is generated by the sets
 $S_0(R)$ and} $S_0(R)^{-1}:= \{ s^{-1} \, | \, s\in S_0(R)\}$.
 \item $Q_l(R)^* = \{ s^{-1}t\, | \, s,t\in S_0(R)\}$.
 \item $Q_l(Q_l(R))=Q_l(R)$.
\end{enumerate}

\end{theorem}

{\it Proof}. 1--3. It is obvious that $G:=\langle S_0(R),
S_0(R)^{-1}\rangle\subseteq Q_l(R)^*\subseteq S_0(Q_l(R))$.
Applying Proposition \ref{A4Jul10}.(2) in the situation where
$S=S_0(R)$ and $T = S_0(Q_l(R))$ we see that $$S_0(R)\subseteq
T':= S_0(Q_l(R))\cap R\in \Den_l(R, 0),$$ and so $S_0(R) = T'$, by
the maximality of $S_0(R)$. Let $q\in S_0(Q_l(R))$. Then $q=
s^{-1} t$ for some elements $s\in S_0(R)$ and $t = sq\in
S_0(Q_l(R))\cap R= S_0(R)$. Therefore, $S_0(Q_l(R)) \subseteq \{
s^{-1} t\, | \, s,t\in S_0(R)\} \subseteq G$, and so $G= Q_l(R)^*
= S_0(Q_l(R)) = \{ s^{-1} t\, | \, s,t \in S_0(R)\}$. This proves
statements 1--3.

4. Statement 4 follows from statement 1. $\Box $

$\noindent $

{\bf Necessary and sufficient conditions for  $Q_l(R)$ to be  a
semi-simple ring}. A ring $Q$ is called a {\em ring of quotients}
if every element $c\in \CC_Q$ is invertible. A subring $R$ of a
ring of quotients $Q$ is called a {\em left order} in $Q$ if
$\CC_R$ is a left Ore set and $\CC_R^{-1}R=Q$. A ring $R$ has {\em
finite left rank} (i.e. {\em finite left uniform dimension}) if
there are no infinite direct sums of nonzero left ideals in $R$.

The next theorem gives an answer to the question of when $Q_l(R)$
is a semi-simple ring.

\begin{theorem}\label{5Jul10}
The following properties of a ring $R$ are equivalent.
\begin{enumerate}
\item  $Q_l(R)$ is a semi-simple ring. \item $Q_{cl}(R)$  exists
and is a semi-simple ring. \item $R$ is a left order in a
semi-simple ring. \item $R$ has finite left rank, satisfies the
ascending chain condition on left annihilators and is a semi-prime
ring. \item A left ideal of $R$ is essential iff it contains a
regular element.
\end{enumerate}
If one of the equivalent conditions hold then $S_0(R) = \CC_R$ and
$Q_l(R) = Q_{cl}(R)$.
\end{theorem}

{\it Proof}. Goldie's Theorem states that $2\Leftrightarrow 3
\Leftrightarrow 4\Leftrightarrow 5$.

$(3\Rightarrow 1)$ If statement 3 holds then $Q_l(R) =
\CC_R^{-1}R$ is a semi-simple ring.

$(1\Rightarrow 3)$ We have to show that $\CC_R$ is a left Ore set
since then $\CC_R = S_0(R)$ and $Q_l(R) = Q_{cl}(R)$ is a
semi-simple ring. Notice that $S_0(R)\subseteq \CC_R$. Let $q\in
\CC_R$. Then the $Q_l(R)$-module monomorphism $\cdot q: Q_l(R) \ra
Q_l(R)$, $x\mapsto xq$, is an isomorphism, and its inverse is
necessarily of the form $\cdot p$ for some element $p\in Q_l(R)$.
Clearly, $pq=1$ and $qp=1$, i.e.  $q\in Q_l(R)^*$. By Theorem
\ref{4Jul10}.(3), $q=s^{-1}t$ for some elements $s,t\in S_0(R)$.
We have to show that, for each $q\in \CC_R$ and  $r\in R$, there
exists an element $c\in \CC_R$ such that $crq^{-1} = crt^{-1}s\in
R$. Since $t\in S_0(R)$, there exists an element $s_1\in S_0(R)$
such that $s_1rt^{-1} \in R$. It suffices to take $c=s_1$ since
$S_0(R)\subseteq \CC_R$.  $\Box $

$\noindent $

The next corollary gives an interesting criterion of when the
classical quotient ring $Q_{cl}(R)=\CC_R^{-1}R$ exists.

\begin{corollary}\label{A5Jul10}
If the ring $Q_l(R)$ is a left artinian ring then $S_0(R)= \CC_R$
and $Q_{cl}(R) = Q_l(R)$.
\end{corollary}

{\it Proof}.  Each left artinian ring is left Noetherian, hence
the left $Q_l(R)$-module $Q_l(R)$ has finite length. Now, we can
repeat the proof of the implication $1\Rightarrow 3$ of  Theorem
\ref{5Jul10}, where, in fact, we {\em only} used the fact that the
left $Q_l(R)$-module $Q_l(R)$ has finite length. $\Box $


\begin{proposition}\label{Pr11.6Sten}
{\rm (Proposition 11.6, \cite{Stenstrom-RingsQuot};
\cite{Jondrup})} Let $A$ be a subring of a ring $B$. If $M$ is a
finitely generated flat $A$-module such that $B\t_AM$ is a
projective $B$-module then $M$ is a projective $A$-module.
\end{proposition}

\begin{corollary}\label{11Jul10}
If there exists a finitely generated flat $R$-module $M$ which is
not projective the the ring $Q_l(R)$ is not a semi-simple ring.
\end{corollary}

{\it Proof}. If $Q_l(R)$ were a semi-simple ring then
$Q_l(R)\t_RM$ would be a projective $Q_l(R)$-module, and so $M$
would be a projective $R$-module, by Proposition \ref{Pr11.6Sten},
a contradiction. $\Box $



\section{The maximal left quotient rings of
a ring}\label{TMLQR}

In this section, a new class of rings, the class of left
localization maximal rings, is introduced. It is proved that, for
an arbitrary ring $R$, the set of maximal elements of the poset
$(\Den_l(R), \subseteq )$ is a non-empty set (Lemma
\ref{b14Nov10}.(2)), and therefore the set of maximal left
quotient rings of the ring $R$ is a non-empty set. A criterion is
given (Theorem \ref{21Nov10}) for a left quotient ring of a ring
to be a maximal left quotient ring of the ring. Many results on
denominator sets are proved. In particular, for each denominator
set $S\in \Den_l(R, \ga )$, connections are established between
the left denominator sets $\Den_l( R, \ga )$, $\Den_l(R/ \ga , 0)$
and $\Den_l(S^{-1}R, 0)$.

\begin{proposition}\label{c3Jul10}
\begin{enumerate}
\item For each ring $A\in \Loc_l(R, \ga )$ where $\ga \in
\Ass_l(R)$, the set $\Den_l(R, \ga , A) := \{ S\in \Den_l(R, \ga )
\, | \, S^{-1}R=A\}$ is an ordered submonoid of $\Den_l(R, \ga )$,
and  \item $S(R, \ga , A) := \bigcup_{S\in \Den_l(R, \ga , A)}S$
is its largest element. In particular, $S_0(R) = S(R, 0, Q_l(R))$.
\item Let $S_i\in \Den_l(R, \ga , A)$, $i\in I$, where $I$ is an
arbitrary non-empty set. Then $\langle S_i\, | \, i\in I\rangle\in
\Den_l(R, \ga , A)$ (see (\ref{SiIgen})). Moreover, $\langle S_i\,
| \, i\in I\rangle$  is the least upper bound of the set $\{
S_i\}_{i\in I}$ in $\Den_l(R, \ga , A)$ and in $\Den_l(R, \ga )$.
\end{enumerate}
\end{proposition}

{\it Proof}. 1. In view of Theorem \ref{3Jul10}, it suffices to
show that if $S_1, S_2\in \Den_l(R, \ga , A)$ then $S_1S_2\in
\Den_l(R, \ga , A)$, i.e. $(S_1S_2)^{-1}R=A$. By Theorem
\ref{3Jul10}, $S_1S_2\in \Den_l(R, \ga )$. Notice that $A=
S_1^{-1} R=S_2^{-1} R$. For each $s=s_1\cdots s_n\in S_1S_2$ where
all $s_i\in S_1\cup S_2$, and, for each $r\in R$,
$$ s^{-1} r= s_n^{-1} \cdots s_1^{-1} r\in s_n^{-1}\cdots
(s_2^{-1}A) \subseteq s_n^{-1}\cdots (s_3^{-1}A) \subseteq \cdots
\subseteq A.$$ Therefore, $(S_1S_2)^{-1} R = A$.

2. By Theorem \ref{3Jul10}, $S(R, \ga , A) = \langle S\, | \, S\in
\Den_l(R, \ga , A)\rangle\in \Den_l(R,\ga )$; and $S(R, \ga ,
A)^{-1}R=\injlim S^{-1}R= \injlim A= A$ where the injective limit
 is over $S\in \Den_l(R, \ga , A)$.

3. Repeat the proof of statement 2 replacing $\Den_l(R, \ga , A)$
by $I$.  $\Box $


\begin{lemma}\label{a5Jul10}

\begin{enumerate}
\item Let $S\in \Den_l(R, \ga )$, $\gb$ be an ideal of the ring
$R$ such $\gb \subseteq \ga $, and $\pi : R\ra R/ \gb$, $a\mapsto
\oa = a+\gb$. Then $\pi (S) \in \Den_l (R/ \gb , \ga / \gb )$ and
$S^{-1} R\simeq \pi (S)^{-1} (R/ \gb )$. \item Let $S_1, S_2\in
\Den_l(R)$ and $S_1\subseteq S_2$. Then
\begin{enumerate}
\item $\ga_1:= \ass (S_1) \subseteq \ga_2:= \ass (S_2)$; there is
the $R$-ring homomorphism $\v : S_1^{-1} R\ra S_2^{-1}R$, $s^{-1}
a\mapsto s^{-1} a$; and $\ker (\v ) = S_1^{-1} (\ga_2/ \ga_1)$.
\item Let $\pi_1 : R\ra R/ \ga_1$, $a\mapsto \oa = a+\ga_1$, and
$\ttS_2$ be the multiplicative submonoid of $(S_1^{-1}(R/\ga _1),
\cdot )$ generated by $\pi_1(S_2)$ and $\pi_1(S_1)^{-1}= \{
\bs^{-1} \, | \, s\in S_1\}$. Then $\pi_1(S_2), \ttS_2\in
\Den_l(S_1^{-1}R, S_1^{-1} (\ga_2/ \ga_1))$ and $\ttS_2^{-1}
(S_1^{-1}R) \simeq \pi_1(S_2)^{-1} (S_1^{-1}R)\simeq S_2^{-1} R$.
\end{enumerate}
\end{enumerate}
\end{lemma}

{\it Proof}. 1. Since $\pi$ is an epimorphism and $ \gb \subseteq
\ga$, $\pi (S)$ is a left Ore set in $R/ \ga$. Clearly, $\ga /
\gb\subseteq \ass (\pi (S))$. To prove that the inverse inclusion
holds, let $\oa\in \ass (\pi (S))$. Then $0= \bs \oa$ for some
element $\bs \in \pi (S)$. Then $b:= sa\in \gb$. Since
$\gb\subseteq \ga$, we can find an element $t\in S$ such that
$tb=0$. Then $(ts) a=0$, and so $a\in \ga$. Therefore, $\ass (\pi
(S)) = \ga / \gb$. To prove that $\pi (S) \in \Den_l(R/ \gb )$, it
remains us to show that $\oa \bs =0$, for some $\oa \in R/ \gb$
and $\bs \in \pi (S)$, implies that $\ot \oa =0$ for some $t\in
S$. Clearly, $b:=as \in \gb \subseteq \ga$, and so $0=
s_1b=(s_1a)s$ for some $s_1\in S$. It follows that $s_1 a\in \ga$
since $S\in \Den_l(R, \ga )$. We can find an element $s_2\in S$
such that $s_2s_1a=0$. Therefore, $\pi (S)\in \Den_l(R/ \gb, \ga /
\gb )$. It suffices to take $t= s_2s_1\in S$.

By the universal property of the ring $S^{-1}R$, there is the ring
epimorphism $$S^{-1}R\ra \pi (S)^{-1}(R/ \gb ), \;\; s^{-1}
r\mapsto \pi (s)^{-1}\pi (r).$$  By the universal property of the
ring $\pi (S)^{-1}(R/ \gb )$, there is the ring epimorphism $ \pi
(S)^{-1}(R/ \gb )\ra S^{-1}R$, $ \pi (s)^{-1}\pi (r)\mapsto s^{-1}
r$. Therefore, $S^{-1}R\simeq \pi (S)^{-1}(R/ \gb )$.

2a. The inclusion $\ga_1\subseteq \ga_2$ is obvious. The image
$\bS_1$  of the Ore set $S_1$ under the ring epimorphism $\pi :
R\ra R/ \ga_2$ is an Ore set in $R/ \ga_2$ since $S_1\subseteq
S_2$ and $S_2\in \Den_l( R, \ga_2)$. Therefore, each element of
the subring of $S_2^{-1}R$ generated by $\pi (R)$ and $\bS_1$ has
the form $\pi (s)^{-1}\pi (r)$ for some $s\in S_1$ and $r\in R$.
By the universal property of the ring $S_1^{-1}R$, the map $\v $
exists. An element $s^{-1}r\in S_1^{-1}R$, where $s\in S_1$ and
$r\in R$, belongs to the kernel of $\v$ iff $s^{-1} r=0$ in
$S_2^{-1}R$ iff $tr=0$ for some $t\in S_2$ iff $r\in \ga_2$.
Therefore, $\ker (\v ) = S_1^{-1} (\ga_2/ \ga_1 )$.

2b. First, we prove that the set $\pi_1(S_2)$ (resp. $\ttS_2$) is
a left Ore set in $S_1^{-1}R$, i.e. we have to prove that  for
each element $\bs_2\in \pi_1(S_2)$ (resp. $\bs_2\in \ttS_2$) and
$s_1^{-1} r \in S_1^{-1}R$ where $s_1\in S_1$ and $r\in R$, there
exist elements $t\in \pi_1(S_2)$ (resp. $t\in \ttS_2$) and $a\in
S_1^{-1}R$ such that $ts_1^{-1} r= a\bs_2$. In the ring $S_2^{-1}
R$, the  element $s_1^{-1} rs_2^{-1}$ can be written as $ s_3^{-1}
r_1$ for some $s_3\in S_2$ and $r_1\in R$. In $S_1^{-1}R$, $b:=
\bs_3s_1^{-1} r-r_1\bs_2\in \ker (\v ) = S_1^{-1} (\ga_2/ \ga_1)$,
by statement 2(a). Therefore, there exists an element $s_4\in S_2$
such that $\bs_4b =0$, and so $\overline{s_4s_3}\cdot s_1^{-1} r =
\bs_4r\cdot \bs_2$. It suffices to take $t=\overline{s_4s_3}\in
\pi_1(S_2) \subseteq \ttS_2$ and $a= \bs_4r$.

It is obvious that $S_1^{-1} (\ga_2/ \ga_1)\subseteq \ass
(\pi_1(S_2))\subseteq \ass (\ttS_2)$. If $u\in \ass (\ttS_2)$ then
$su=0$ for some $s\in \ttS_2$. Then $0=\v (su)= \v (s) \v (u)$,
and so $u\in \ker (\v ) = S_1^{-1} (\ga_2/ \ga_1)$ (by statement
2(a)) since $\v (s)$ is a unit. Therefore, $S_1^{-1} (\ga_2/
\ga_1)=\ass (\pi_1(S_2))=\ass (\ttS_2)$.

If $vs=0$ for some elements $v\in S_1^{-1}R$ and $s\in \pi_1(S_2)$
(resp. $s\in \ttS_2$). Then $0=\v (v) \v (s)$, and so $\v (v)\in
\ker (\v )= S_1^{-1} (\ga_2/ \ga_1)$ since $\v (s)$ is a unit.
This proves  that $\pi_1(S_2) \in \Den_l(S_1^{-1}R, S^{-1}(\ga_2/
\ga_1 ))$ (resp. $\ttS_2\in \Den_l (S_1^{-1}R, S^{-1}(\ga_2/ \ga_1
)))$.  Now, it is obvious that $\pi_1(S_2)^{-1} (S_1^{-1} R)
\simeq S_2^{-1}R \simeq \ttS_2^{-1} (S_1^{-1}R)$. $\Box $

$\noindent $

The set $(\Loc_l(R, \ga ),  \ra )$ is a poset where $A_1\ra A_2$
if $A_1 = S_1^{-1}R$ and  $A_2= S_2^{-1}R$ for some denominator
sets $S_1, S_2 \in \Den_l(R, \ga )$ with $S_1\subseteq S_2$,  and
$A_1\ra A_2$ is the inclusion map in Lemma \ref{a5Jul10}.(2a). If
$(S_1', S_2')$ is another such a pair then, by Proposition
\ref{c3Jul10}.(1),  $A_1= S_1^{-1}R = S_1'^{-1}R= (S_1S_1')^{-1}
R\ra A_2= S_2^{-1} R = S_2'^{-1}R = (S_2S_2')^{-1}R$; $S_1S_1',
S_2S_2'\in \Den_l(R, \ga )$ with $S_1S_1'\subseteq S_2S_2'$.

$\noindent $

In the same way, the poset $(\Loc_l(R), \ra )$ is defined, i.e.
$A_1\ra A_2$ if there exists $S_1, S_2\in \Den_l(R)$ such that
$S_1\subseteq S_2$, $A_1=S_1^{-1}R$ and $A_2=S_2^{-1}R$, $A_1\ra
A_2$ stands for the map $\v :S_1^{-1}R\ra S_2^{-1}R$ (Lemma
\ref{a5Jul10}.(2a)). The map 
\begin{equation}\label{m1RD}
(\cdot )^{-1} R: \Den_l(R)\ra \Loc_l(R), \;\; S\mapsto S^{-1}R,
\end{equation}
is an epimorphism of the posets $(\Den_l(R), \subseteq )$ and
$(\Loc_l(R), \ra )$. For each ideal $\ga \in \Ass_l(R)$, it
induces the epimorphism of posets $(\Den_l(R, \ga ), \subseteq )$
and $(\Loc_l(R, \ga ), \ra )$, 
\begin{equation}\label{m1RD1}
(\cdot )^{-1} R: \Den_l(R, \ga )\ra \Loc_l(R, \ga ), \;\; S\mapsto
S^{-1}R.
\end{equation}

The sets $\Den_l(R)$ and $\Loc_l(R)$ are the disjoint unions

\begin{equation}\label{Den1}
\Den_l(R)=\bigsqcup_{\ga \in \Ass_l(R)}\Den_l(R, \ga ), \;\;
\Loc_l(R)=\bigsqcup_{\ga \in \Ass_l(R)}\Loc_l(R, \ga ).
\end{equation}
For each ideal $\ga \in \Ass_l(R)$, the set $\Den_l(R, \ga )$ is
the disjoint union 
\begin{equation}\label{Den2}
\Den_l(R, \ga ))=\bigsqcup_{A \in \Loc_l(R, \ga ))}\Den_l(R, \ga ,
A ).
\end{equation}
Let $\LDen_l(R, \ga ) :=\{ S(R, \ga , A) \, | \, A\in \Loc_l(R,
\ga ) \}$, see Proposition \ref{c3Jul10}.(2). The map
\begin{equation}\label{Den3}
(\cdot )^{-1}R : \LDen_l( R, \ga ) \ra \Loc_l(R, \ga ) , \;\;
S\mapsto S^{-1} R,
\end{equation}
is an isomorphism of posets.

%

$\noindent $

For a left denominator set $S\in \Den_l(R, \ga )$, there are
natural ring homomorphisms $$R\ra R/\ga \ra S^{-1}R.$$ Lemma
\ref{a6Jul10} and Proposition \ref{8Jul10} establish connections
between the left denominator sets $\Den_l(R, \ga )$, $\Den_l(R/
\ga , 0)$ and $\Den_l(S^{-1}R, 0)$.

Let $S,T\in \Den_l(R)$. The denominator set $T$ is called $S$-{\em
saturated} if $sr\in T$, for some $s\in S$ and $r\in R$, then
$r\in T$, and if $r's'\in T$, for some $s'\in S$ and $r'\in R$,
then  $r'\in T$.

\begin{lemma}\label{a6Jul10}
Let $S\in \Den_l(R, \ga )$, $\pi : R\ra R/\ga$, $a\mapsto a+\ga$,
and  $\s : R\ra S^{-1}R$, $ r\mapsto r/1$.
\begin{enumerate}
\item  Let $T\in \Den_l(S^{-1}R, 0)$ be such that $\pi (S), \pi
(S^{-1})\subseteq T$. Then $T':= \s^{-1} (T) \in \Den_l(R, \ga )$,
$T'$ is $S$-saturated, $T=\{ s^{-1}t'\, | \, s\in S, t'\in T'\}$,
and $S^{-1} R\subseteq T'^{-1}R= T^{-1}R$. \item $\pi^{-1} (S_0(R/
\ga )) = S_\ga (R)$, $ \pi (S_\ga (R)) = S_0(R/\ga ))$ and $ Q_\ga
(R)=S_\ga (R)^{-1}R = Q_l(R/\ga )$.
\end{enumerate}

\end{lemma}

{\it Proof}. 1. In the proof below we often identify the elements
$\s (s)$, $s/1$ and $s$ in order to avoid cumbersome notation. By
the very definition, $T'$ is a multiplicatively closed set that
contains $S$. To prove that $T'$ is a left Ore set in $R$ we have
to show that, for any $t'\in T'$ and $r\in R$, there exist
elements $t_1'\in T'$ and $a\in R$ such that $t_1'r= at'$. Since
$T\in \Den_l(S^{-1}R, 0)$ and $\s (t')\in T$, $t\s (r) = s_1^{-1}
r_1\s (t')$ for some elements $t\in T$, $s_1\in S$ and $r_1\in R$.
The element $t$ is equal to $s_2^{-1} r_2$ for some elements
$s_2\in S$ and $r_2\in R$. Clearly, $r_2\in T'$ since $\pi
(S)\subseteq T$ and $r_2/1 = s_2t\in T$. Fix an element $s_3\in S$
such that $r_3/1:= s_3s_1s_2^{-1} \in \pi (R)$. Then $r_3\in T'$
since $\pi (S)^{\pm 1}\subseteq T$. Multiplying the equality
$$s_1s_2^{-1}r_2\s (r) = r_1\s (t')$$
 by the element $s_3$ on the
left we obtain the equality $r_3r_2\s (r) = s_3r_1\s (t')\in \pi
(R)$. Hence $s_4r_3r_2\cdot r= s_4s_3r_1\cdot t'$ for some element
$s_4\in S$. Now, take $t_1':= s_4r_3r_2\in T'$ and
$a:=s_4s_3r_1\in R$.

 Since
$S\subseteq T'$, $\ga \subseteq \ass (T')$. If $\th u=0$ and
$v\eta =0$ for some elements $u,v\in R$ and $\th , \eta \in T'$
then $0=\s (\th ) \s (u) $ and $0= \s (v)\s (\eta )$ and so $0= \s
(u) = \s (v)$ since the elements $\s (\th ) $ and $\s (\eta )$ are
units. Then, $u,v\in \ga$. This proves that $T'\in \Den_l(R,\ga
)$.

 The denominator set $T'$ is $S$-saturated since $sr = t'\in
T'$ (resp. $rs=t'\in T'$)  for some $s\in S$ and $r\in R$ implies
$r/1= s^{-1} t'\in T$ (resp.  $r/1= t' s^{-1} \in T$), and so
$r\in T'$. Since $\pi (S)\subseteq T$, $T=\{ s^{-1} t'\, | \, s\in
S, t'\in T'\}$. Then  the inclusion and the equality,
$S^{-1}R\subseteq T'^{-1}R = T^{-1}R$, are obvious.

2. Let $S':= \pi^{-1} (S_0 (R/ \ga ))$. Since the map $\pi$ is
surjective we have $\pi (S') = S_0(R/\ga )$. We aim to show that
$S'= S_\ga (R)$.

{\em Step 1}: $S'\in \Ore_l(R)$. We have to show that for any
$s\in S'$ and $r\in R$ there exist elements $s_1\in S'$ and $
r_1\in R$ such that $s_1 r= r_1s$. Since $\pi (S') = S_0(R/ \ga
)\in \Ore_l(R/ \ga )$, we have $\pi (s') \pi ( r) = \pi (r') \pi
(s)$ for some elements $s'\in S'$ and $r'\in R$. Then $\pi ( s'
r-r's)=0$, and so $s'r-r's\in \ga$, hence $s''(s'r-r's)=0$ for
some $s''\in S$. Note that $\pi (S) \in \Den_l(R/ \ga , 0)$ hence
$\pi (S) \subseteq S_0(R/ \ga )$, and so $S\subseteq S'$. Now,
$s_1r=r_1s$ where $s_1= s''s'\in S'$ and $ r_1 = s''r'\in R$.
Therefore, $S'\in \Ore_l(R)$.

{\em Step 2}: $\ass (S') = \ga $. Since $S\subseteq S'$, we have
the inclusion $\ga \subseteq \ass (S')$. Let $ r\in \ass (S')$.
Then $s' r=0$ for some $s'\in S'$, and so $\pi ( s') \pi (r) =0$.
It follows that $\pi (r)=0$ since $\pi ( s') \in S_0( R/\ga )$,
i.e. $r\in \ker (\pi ) = \ga $. This means that $\ass (S') = \ga$.

{\em Step 3}: $S'\in \Den_l(R, \ga )$.  In view of Steps 1 and 2,
 we have to show that the equality  $rs'=0$ for some elements $r\in
R$ and $ s' \in S'$ implies that $s'' r=0$ for some $s''\in S'$.
Since $\pi ( r ) \pi ( s')=0$ and $\pi ( s')\in S_0(R/ \ga ) = \pi
( S')$ we have the equality $ \pi ( r ) =0$ and so $ r\in \ga$.
Then $s''r=0$ for some element $s''\in S\subseteq S'$.

{\em Step 4}: $S'= S_\ga (R)$. By Step 3, $S'^{-1} R = \pi (
S')^{-1} R/\ga = S_0(R/ \ga )^{-1} R/ \ga = Q_l( R/ \ga )$. Notice
that $S_\ga (R)^{-1}R = \pi (S_\ga (R))^{-1} R/ \ga \subseteq
Q_l(R/ \ga )$. Since $\pi (S_\ga (R)) \in \Den _l( R/ \ga , 0)$,
we have $\pi ( S_\ga (R)) \subseteq S_0 (R/ \ga )$, hence $ S_\ga
(R) \subseteq \pi^{-1} (S_0 (R/ \ga )) = S' \subseteq S_\ga (R)$,
i.e. $S' = S_\ga (R)$. $\Box $

$\noindent $

For $S_1, S_2\in \Den_l(R)$ such that $S_1\subseteq S_2$,
$[S_1,S_2]:= \{ T\in \Den_l(R)\, | \, S_1\subseteq T\subseteq
S_2\}$ is an {\em interval} in the posed $\Den_l(R)$. If, in
addition, $S_1, S_2\in \Den_l(R, \ga )$ then $[S_1, S_2]\subseteq
\Den_l(R, \ga )$ since $S_1\subseteq S\subseteq S_2$ implies $\ga
= \ass (S_1) \subseteq \ass (S)\subseteq \ass (S_2)=\ga $, i.e.
$\ass (S) = \ga$.

\begin{proposition}\label{8Jul10}
Let $S\in \Den_l(R, \ga )$; $\pi :R\ra R/ \ga$, $a\ra \oa =
a+\ga$;
 $\s : R\ra S^{-1}R, r\ra r/1$; $G:=\langle \pi (S), \pi
 (S)^{-1}\rangle \subseteq (S^{-1}R)^*$ (i.e. $G$ is the subgroup
 of the group $(S^{-1}R)^*$ of units of the ring $S^{-1}R$
 generated by $\pi (S)^{\pm 1}$).
\begin{enumerate}
\item Let  $ [\s^{-1} (G), S_\ga (R)]_{S-{\rm
 com}}:= \{ S_1\in [\s^{-1} (G), S_\ga (R)]\, | \, \s^{-1} (G\pi (S_1))=S_1\}$ and
 $[ G, S_0(S^{-1}R)]:= \{ T \in \Den_l(S^{-1}R, 0)\, | \, G\subseteq T\subseteq S_0(S^{-1}R)\}$. Then the map
$$[\s^{-1}(G), S_\ga (R)]_{S-{\rm com}}\ra [ G, S_0(S^{-1}R)], \;\;
S_1\mapsto \ttS_1:= G\pi (S_1), $$ is an isomorphism of posets and
abelian monoids with the inverse map $T\mapsto \s^{-1}(T)$ where
$G\pi (S_1)$ is the multiplicative monoid generated by $G$ and
$\pi (S_1)$ in $S^{-1}R$. In particular,
$$G\pi (S_\ga (R)) = S_0(S^{-1}R), \;\; S_\ga(R) = \s^{-1}
(S_0(S^{-1}R)), \;\; S_\ga (R)^{-1}R = Q_l(R/\ga ),$$ the monoid
$[\s^{-1}(G), S_\ga (R)]_{S-{\rm
 com}}$ is a complete lattice  (since $[ G, S_0(S^{-1}R)]$ is a complete  lattice, as an interval
 of the complete lattice $\Den_l(S^{-1}R, 0)$, Corollary \ref{a7Jul10}), and the
 map $S_1\mapsto \widetilde{S}_1$ is a lattice isomorphism.
 \item Consider the interval $[G\cap (R/ \ga ) , S_0(R/ \ga )]$ in
 $\Den_l(R/ \ga  , 0)$. Let $[G\cap (R/ \ga ) , S_0(R/ \ga
 )]_{G-{\rm com}}:=\{ T\in [G\cap (R/ \ga ) , S_0(R/ \ga )]\, | \,
 GT\cap (R/ \ga ) =T\}$. Then $[G\cap (R/ \ga ) , S_0(R/\ga )]_{G-{\rm com}} \subseteq
 \Den_l(S^{-1}R, 0)$ and  the map
 $$[G\cap (R/ \ga ) , S_0(R/ \ga
 )]_{G-{\rm com}}\ra [G, S_0(S^{-1}R)], \;\; T\mapsto GT,$$
 is an isomorphism of posets and abelian monoids with the inverse
 map $T'\mapsto T'\cap (R/ \ga )$ where $GT$ is the product in $\Den_l(S^{-1}R, 0)$. In particular,
 the monoid $[G\cap (R/ \ga ) , S_0(R/ \ga
 )]_{G-{\rm com}}$ is a complete lattice.
  \item The map $$[ \s^{-1} (G), S_\ga (R)]_{S-{\rm com}}\ra [G\cap (R/ \ga ) , S_0(R/ \ga
 )]_{G-{\rm com}}, \;\; S_1\mapsto G\pi (S_1)\cap (R/ \ga ),$$
is an isomorphism of posets and abelian monoids with the inverse
 map $S'\mapsto\pi^{-1}(S')$.
\end{enumerate}
\end{proposition}

{\it Proof}. 1. The equality $\widetilde{\s^{-1}(G)}=G$ is
obvious. Then, by Lemma \ref{a5Jul10}.(2).(b), the map $\phi
:S_1\mapsto \ttS_1$ is well-defined. By Lemma \ref{a6Jul10}.(1),
the map $\psi : T\mapsto \s^{-1} (T)$ is well-defined and
$\widetilde{\s^{-1}(T)}=T$ (i.e. $\phi \psi = 1$) since $T=\{
s^{-1}t\, | \, s\in S, t\in \s^{-1}(T)\}$ (Lemma
\ref{a6Jul10}.(1)) and $G\subseteq T$. By the very definition of
the set $[S, S_\ga (R)]_{S-{\rm com}}$, $\psi \phi =1$, i.e. $\psi
= \phi^{-1}$.
 For all elements $S_1, S_2\in [S, S_\ga
(R)]_{S-{\rm com}}$,
$$ \widetilde{S_1\bigwedge S_2}= \widetilde{S_1
S_2}= G\pi (S_1)\pi (S_2)=  G\pi (S_1)G\pi (S_2) = \ttS_1\ttS_2 =
\ttS_1\bigwedge \ttS_2,$$ and so the map $\phi$ is an isomorphism
of posets and abelian monoids. Therefore, $G\pi (S_\ga (R))
=S_0(S^{-1}R)$ and $S_\ga (R) = \s^{-1} (S_0(S^{-1}R))$. Then,
$S_\ga (R)^{-1}R = Q_l(R/ \ga )$, by Lemma \ref{a6Jul10}.(1).  The
rest is obvious.

2. Let $\phi : T\mapsto GT$ and $ \psi :T'\mapsto T'\cap \bR$
where $\bR:= R/\ga $.

{\em Step 1: $\phi$ is well-defined}. Since $G\in \Den_l(S^{-1}R ,
0)$ and $\phi (T) = GT$ (the product in $\Den_l(S^{-1}R, 0)$), we
have to show that $T\in \Den_l(S^{-1}R , 0)$.

First, let us show that $T\in \Ore_l(S^{-1}R)$, i.e. for each
$s^{-1}r\in S^{-1}R$ (where $s\in S$ and $r\in \bR$) and $t\in T$
we have to show that $t_1s^{-1} r= at$ for some elements $t_1\in
T$ and $ a\in S^{-1}R$. Since $T\in  \Den_l(\bR , 0)$, $ t'r=at$
for some $t'\in T$ and $a\in \bR$. So, it suffices to take
$t_1=t's$ ($ t_1\in T$ since $\pi (S) \subseteq T$).

Let us show that $\ass (T)=0$ in $S^{-1}R$. Suppose that $t\cdot
s^{-1} r=0$ for some $t\in T$, $s\in S$ and $r\in R$, we have to
show that $s^{-1} r=0$. There exist elements $s_1\in S$ and
$t_2\in \bR$ such that $s_1t= t_2s$, hence $t_2=s_1ts^{-1}\in \bR
\cap GT=T$. Then $0=s_1\cdot 0 = s_1ts^{-1} r= t_2r$, hence $r=0$
(since $T\in \Den_l( \bR , 0)$), and so $ s^{-1}r=0$, as required.
Therefore, $\ass (T)=0$ in $S^{-1}R$. To finish the proof of Step
1, we have to prove that $s^{-1} rt=0$ for some $s\in S$, $r\in R$
and $ t\in T$ implies $s^{-1}r=0$. This is obvious since $s^{-1}
rt=0$ implies $rt=0$ in $\bR$, hence $r=0$ since $T\in \Den_l(\bR
, 0)$, and so $s^{-1} r=0$. Therefore, $T\in \Den_l(S^{-1}R, 0)$
for all $T\in [ G\cap (R/ \ga ), S_0(R/\ga )]_{G-{\rm com}}$. In
particular, $[ G\cap (R/ \ga ), S_0(R/\ga )]_{G-{\rm
com}}\subseteq \Den_l(S^{-1}R, 0)$.

{\em Step 2: $\psi$ is well-defined}. Let $T'\in \Den_l(S^{-1}R,
0)$. We have to show that $T:=\bR \cap T'\in \Den_l(\bR , 0)$. It
is obvious that $ \ass_{\bR } (T)=0$ since $\ass_{\bR }(T)
\subseteq \ass_{S^{-1} R} (T')=0$. If $rt=0$ for some elements
$r\in \bR$ and $t\in T$ then  $r=0$ since $T\subseteq T'$, $T' \in
\Den_l(S^{-1}R, 0)$ and $\bR \subseteq S^{-1}R$. It remains to
show that $T\in \Ore_l(\bR )$, that is, for any $r\in \bR$ and
$t\in T$,  $t_1r=r_1t$ for some elements $r_1\in \bR$ and $ t_1\in
T$.  Since $T\subseteq T'$ and $T'\in \Den_l(S^{-1}R, 0)$,
$t'r=s_1^{-1}r_1't$ for some elements $t'\in T'$, $s_1\in S$ and $
r_1'\in R$. Fix an element $s_2\in S$ such that $t_1:= s_2s_1t'\in
\bR$. Then $t_1\in \bR \cap T'= T$ and $t_1r=s_2s_1 t'r=
s_2s_1s_1^{-1} r_1't = s_2r_1'\cdot t= r_1 t$ where $r_1:= s_2
r'\in \bR$, as required.

 By the very definition of the set $[
G\cap \bR , S_0(\bR )]$, $\psi \phi =1$. To finish the proof of
statement 2 it remains to show that $\phi \psi =1$, i.e. $G(\bR
\cap T') = T'$ for all $T'\in [ G, S_0(S^{-1}R)]$. Since $G, \bR
\cap T' \subseteq T'$, the inclusion $T'':= G(\bR \cap T')
\subseteq T'$ is obvious. The reverse inclusion follows from the
fact that any element $t'\in T'$ can be written as $s^{-1} t$ for
some elements $s\in S\subseteq G$ and $t\in \bR$, hence $t=st'\in
\bR \cap T'=T$.

3. Statement 3 follows from statements 1 and 2.  $\Box $

$\noindent $

The elements of the set $[S, S_\ga (R)]_{S-{\rm com}}$ are called
$S$-{\em complete} and the elements of the set $[G\cap (R/ \ga ) ,
S_0(R/ \ga )]_{G-{\rm com}}$ are called $G$-{\em complete}.

\begin{lemma}\label{a21Nov9}
We keep the notation of Proposition \ref{8Jul10}.(1). If $S_1\in [
\s^{-1} (G), S_\ga (R)]_{S-{\rm com}}$ then $S_1$ is
$S$-saturated.
\end{lemma}

{\it Proof}. Notice that $\s^{-1} (G) \subseteq S_1\subseteq S_\ga
(R)$ and $S_1=\s^{-1} (G\pi ( S_1))$. If $ s_1:= sr \in S_1$
(resp. $s_1:= rs\in S_1$) for some elements $s\in S$ and $r\in R$
then $r/1=s^{-1} s_1\in G\pi ( S_1)$ (resp. $ r/1= s_1s^{-1} \in
G\pi ( S_1)$) hence $r\in \s^{-1} (G\pi ( S_1))= S_1$, i.e. $S_1$
is $S$-saturated. $\Box $

$\noindent $

{\bf The maximal left quotient rings of a ring}.
\begin{lemma}\label{a14Nov10}
Let $S_1\subseteq S_2\subseteq \cdots \subseteq S_i\subseteq
\cdots $  be an ascending chain in $\Den_l(R)$, $\ga_i:= \ass
(S_i)$, $S:=\bigcup_{i\geq 1}S_i$. Then $\ga_1\subseteq
\ga_2\subseteq \cdots \subseteq \ga_i\subseteq \cdots$ is the
ascending chain in $\Ass_l(R)$, $S\in \Den_l(R, \ga )$ where $\ga
:= \bigcup_{i\geq 1} \ga_i$, $S^{-1}R= \injlim S_i^{-1}R$ where
$S_1^{-1}R\ra S_2^{-1}R\ra \cdots \ra S_i^{-1}R\ra \cdots$ (Lemma
\ref{a5Jul10}.(2a)).
\end{lemma}

{\it Proof}. By Lemma \ref{a5Jul10}.(2a), $S\in \Den_l(R, \ga )$.
 For each number $i=1, 2, \ldots $, define the ring homomorphism
 $\phi_i : S_i^{-1} R\ra S^{-1}R$, $s^{-1} r\mapsto s^{-1} r$, and
 $\nu_i : S_i^{-1} R \ra S_{i+1}^{-1}R$, $s^{-1}r\mapsto s^{-1}
 r$. Then $\phi_i = \phi_{i+1} \nu_i$  for all $i$.
Hence, there is the ring homomorphism $\phi : \injlim S_i^{-1}R\ra
S^{-1}R$ which is a surjection since $S= \bigcup_{i\geq 1} S_i$
and has kernel $0$ since $\ga = \bigcup_{i\geq 1} \ga_i$, i.e.
$\phi$ is an isomorphism. $\Box $

$\noindent $

Consider the poset $(\Den_l(R), \subseteq )$. For each element
$S\in \Den_l(R)$, let $[S, \cdot ) := \{ S'\in \Den_l(R)\, | \,
S\subseteq S'\}$.

\begin{lemma}\label{b14Nov10}
\begin{enumerate}
\item  For each element $S\in \Den_l(R)$, there exists a maximal
element in the poset $([S, \cdot ), \subseteq )$. \item The set
$(\maxDen_l(R), \subseteq )$ of maximal elements of the poset
$(\Den_l(R), \subseteq )$ is a non-empty set.
\end{enumerate}
\end{lemma}

{\it Proof}. 1. Statement 1 follows from Lemma \ref{a14Nov10} and
Zorn's Lemma.

2. Statement 2 follows from statement 1 and the fact that the set
$\maxDen_l(R)$ is the set of maximal elements of the poset $[\{
1\} ,\cdot ) = \Den_l(R)$. $\Box $

$\noindent $

{\it Definition}. An element $S$ of the set $\maxDen_l(R)$ is
called a {\em maximal left denominator set} of $R$ and the ring
$S^{-1}R$ is called a {\em maximal left quotient ring} of $R$. The
intersection
$$ \llrad (R) := \bigcap_{S\in \maxDen_l(R)} \ass (S)$$
is called the {\em left localization radical } of the ring $R$.

$\noindent $

\begin{proposition}\label{A21Nov10}
Let $\ga \in \Ass_l(R)$, $Q:= Q_\ga (R)$, $Q^*$ be the group of
units of the ring $Q$ and  $\s : R\ra Q_\ga (R)$, $ r\mapsto
\frac{r}{1}$. Let $T\in \Den_l(Q, \gb )$ where $\gb \in \Ass_l(Q)$
and $Q^*T$ be the multiplicative sub-semigroup of $(Q, \cdot )$
generated by $Q^*$ and $T$. Then
\begin{enumerate}
\item $Q^*T\in \Den_l(Q, \gb )$. \item If, in addition, $Q^*
\subseteq T$ (eg, $Q^*T$ from statement 1) then
\begin{enumerate}
\item $T':= \s^{-1} (T) \in \Den_l(R, \gb')$ where $\gb':= \s^{-1}
(\gb )\supseteq \ga$, $S_\ga (R) \subseteq T'$, $T= Q^* \s (T')$
(i.e. the monoid $T$ is generated by $Q^*$ and $\s (T')$) and
$T'^{-1}R = T^{-1}Q$ (i.e. the natural ring monomorphism
$T'^{-1}R\ra T^{-1}Q$, $ t^{-1}r\mapsto t^{-1}r$, is an
isomorphism). \item $S_\ga (R) \subseteq S_{\gb'}(R)$ and
$S_{\gb'} (R) = \s^{-1} (S_\gb (Q))$. \item $Q_{\gb'}(R) = Q_l( Q/
\gb )$, i.e. the natural ring monomorphism $Q_{\gb'}(R)\ra Q_l( Q/
\gb )$, $ t^{-1} r \mapsto t^{-1}r$, is an isomorphism.
\end{enumerate}
\end{enumerate}
\end{proposition}

{\it Proof}. 1. Clearly, $\CT := Q^* T$ is a multiplicative
monoid. We have to show that $\CT \in \Den_l( Q, \gb )$. An
element $s$ of $\CT$ is a product $s= q_1t_1q_2t_2\cdots q_nt_n$
for some elements $q_i\in Q$ and $t_i\in T$.

 {\it Step 1:  $\CT
\in \Ore_l(R)$}. For any element $s\in \CT$ and $a\in Q$, we have
to find elements $s_1\in \CT$ and $a_1\in Q$ such that $s_1 a=
a_1s$. We use induction on $n$. When $n=1$, i.e. $s= q_1t_1$, then
$s_1a= a't_1$ for some elements $s_1\in T$ and $a'\in Q$ since $
T\in \Den_l( Q, \gb )$. It suffices to take $a_1= a' q_1^{-1}$
since $s_1a= (a'q_1^{-1})\cdot q_1t_1= a_1q_1t_1$. Suppose that
$n>1$ and the statement is true for all $n'<n$. Then $s_n a= a_n
q_n t_n$ for some elements $s_n\in T$ and $a_n\in Q$. By
induction, $s_n'a_n = a_1q_1t_1\cdots q_{n-1}t_{n-1}$ for some
elements $s_n'\in \CT$ and $ a_1\in Q$. It suffices to take $s_1 =
s_n' s_n$ since then
$$ s_1a= s_n' s_n a= s_n' a_n q_nt_n = a_1q_1t_1\cdots
q_{n-1}t_{n-1}q_nt_n= a_1s.$$ {\it Step 2: $\ass (\CT ) = \gb $}.
We have to show that  $sa=0$ for some  $a\in Q$ implies $a\in
\gb$. We use induction on $n$. When $n=1$, i.e. $s= q_1t_1$ then
$q_1t_1 a=0$ implies $t_1 a=0$ (since  $q_1$ is a unit), and so
$a\in \gb$ since $T\in \Den_l( Q, \gb )$. Let $n>1$. Suppose that
the result is true for all $n'<n$. Then $0= sa= q_1t_1\cdots
q_{n-1} t_{n-1} \cdot ( q_nt_n a)$ implies (by induction) $q_n t_n
a\in \gb$, hence $t_na\in \gb$ since $q_n$ is a unit. Then
$t_{n+1}t_n a=0$ for some element $t_{n+1}\in T$ since $ T\in
\Den_l( Q, \gb )$, and finally $a\in \gb$ since $t_{n+1}t_n\in T$.

{\it Step 3: $\CT \in \Den_l( Q, \gb )$}. We have to show that
$as=0$ for some $ a\in Q$ implies $ a\in \gb$. We use induction on
$n$. When $n=1$, i.e. $aq_1t_1=0$, then $aq_1\in \gb$ since $T\in
\Den_l( Q, \gb )$, and so $a\in \gb$ since $q_1$ is a unit. Let
$n>1$. Suppose that the result is true for all $n'<n$. Then $0=as=
(aq_1t_1)\cdot ( q_2t_2\cdots q_nt_n)$ implies (by induction) $
a_1q_1t_1\in \gb$, hence $t_{n+1} aq_1\in \gb$ for some element $
t_{n+1} \in T$ (since $T\in \Den_l( Q, \gb )$), i.e. $t_{n+1} a\in
\gb$ since $ q_1$ is a unit. Therefore, there exists an element
$t_{n+2} \in T$ such that $t_{n+2} t_{n+1} a=0$. This means that
$a\in \gb$ since $t_{n+2}t_{n+1} \in T \in \Den_l( Q, \gb )$.

2a. By the very definition, $T'$ is a monoid and $\ga \subseteq
\gb'$. By Lemma \ref{a6Jul10}.(2), $Q_\ga (R) = Q_0( R/ \ga )$,
$\s^{-1} (S_0(R/ \ga )) = S_\ga (R)$ and $\s ( S_\ga (R)) = S_0(
R/ \ga )$. Therefore, $S_\ga (R) \subseteq T'$ since $\s (S_\ga
(R)) = S_0( R/ \ga )\subseteq Q^*\subseteq T$. Since $Q^*
\subseteq T$, we have the equality $T= Q^* \s (T')$. Let us show
that $T' \in \Den_l( R , \gb ')$.

{\it Step 1: $T' \in \Ore_l( R)$}. We have to show that for any
elements $r\in R$ and $t'\in T'$ there exists elements $t_1'\in
T'$ and $r_1\in R$ such that $t_1'r= r_1t'$. Since $T\in \Den_l(
Q, \gb )$ and $ \s (t')\in T$, we have the equality $t\s (r) =
s^{-1} r_2 \s (t')$ in the ring $Q$ for some elements $t\in T$, $
s\in S_0( R/ \ga )$ and $r_2\in R$. Since $S_0( R/ \ga )\subseteq
Q^*$ and $Q^* \subseteq T$ (by the assumption), the product $st
\in T$. Fix an element $s_1\in S_0( R/ \ga )\subseteq Q^*
\subseteq T$ such that $s_1st=\s (t_2')$ for some element $t_2'\in
R$, necessarily $t_2'\in T'$. Then $s_1 = \s (s_1')$ for some
element $s_1'\in S_\ga (R)$ since $\s (S_\ga (R)) = S_0( R/ \ga
)$. By multiplying the equality $t\s (r) = s^{-1} r_2 \s (t')$ by
the element $s_1s$ on the left, we obtain the equality
$$ \s ( t_2'r-s_1' r_2 t') =0,$$
i.e. $\alpha := t_2'r-s_1'r_2t'\in \ga$. There exists an element
$s_3\in S_\ga (R) \subseteq T'$ such that $s_3\alpha =0$, i.e.
$s_3 t_2'\cdot r= s_3s_1' r_2\cdot t'$. Now, it suffices to take
$t_1':= s_3t_2'\in T'$ and $ r_1 := s_3s_1'r_2\in R$.

{\it Step 2: $\ass (T') = \gb'$}. Let $r\in \ass (T')$, i.e. $t'
r=0$ for some $t'\in T'$. Then $\s (t') \s (r) =0$ in $Q$, and so
$\s (r) \in \gb$ since $\s (t') \in T$. Hence $r\in \gb'$. This
means that $\ass (T') = \gb'$.

{\it Step 3: $T' \in \Den_l( R, \gb ' )$}. We have to show that if
$ rt'=0$ for some $r\in R$ and $t'\in T'$ then $r\in \gb'$. Since
$\s (r) \s (t') =0$ and $ \s (t')\in T$, we have the inclusion $\s
(r)\in \gb$, hence $r\in \gb'$.

Since $\s (T') \subseteq T$ and $\ass (T') = \gb'$, there is the
(natural) ring monomorphism
$$ T'^{-1}R\ra T^{-1}R, \;\; t'^{-1}r\mapsto t'^{-1} r,$$ where
$t'\in T'$ and $r\in R$. The monomorphism is, in fact, an
isomorphism since $T= Q^* \s (T')$, $S_\ga (R)\subseteq T'$ and
$Q^* = \{ \s (s)^{-1} \s (t) \, | \, $ where $s,t\in S_\ga (R)\}$
(by Lemma \ref{a6Jul10}.(2) and Theorem \ref{4Jul10}(1-3)).

2b. Let $T$ and $T'$  be as in statement 2a (it exists, by
statement 1). By statement 2a, $S_\ga (R) \subseteq T'$ and $T'\in
\Den_l( R, \gb')$. Clearly, $T' \subseteq S_{\gb'}(R)$ since $
S_{\gb'} (R)$ is the largest element of the poset $( \Den_l( R,
\gb'), \subseteq )$. Therefore,  $S_\ga (R) \subseteq
S_{\gb'}(R)$.

Let $T= S_\gb (Q)$. Then $Q^*\subseteq T$ by statement 1 and the
maximality of $S_\gb (Q)$. By statement 2a, $T' := \s^{-1} (S_\gb
(Q))\subseteq S_{\gb'} (R)$. Consider the ring epimorphism
$\pi_\ga :R\ra R/ \ga $, $ a\mapsto a+\ga$. Notice that $\ga
\subseteq \gb'$. Applying Lemma \ref{a5Jul10}.(2) in  the
situation $S_\ga (R) \subseteq S_{\gb'}(R)$, we see that $\pi_\ga
(S_{\gb'}(R))\in \Den_l( Q, \gb )$. Therefore,
$S_{\gb'}(R)\subseteq T'$. This proves that $S_{\gb'}(R) = T'$.

2c. Applying statement 2a in the situation $T = S_\gb (Q)
\supseteq Q^*$ and using the fact that $S_{\gb'}(R):= \s^{-1}
(S_\gb (Q)) = \s^{-1} (T)$, we see that $S_{\gb'}(R)^{-1}R= S_{\gb
}(Q)^{-1}Q$, i.e. $Q_{\gb'}(R) = Q_l( Q/ \gb )$ since $Q_l( Q/ \gb
)= S_\gb (Q)^{-1}Q$ (by Lemma \ref{a6Jul10}.(2)). $\Box $

$\noindent $

Let $\maxAss_l(R)$ be the set of maximal elements of the poset
$(\Ass_l(R), \subseteq )$. It is a subset of the set

\begin{equation}\label{mADen}
\assmaxDen_l(R) := \{ \ass (S) \, | \, S\in \maxDen_l(R) \}
\end{equation}
which  is a non-empty set, by Lemma \ref{b14Nov10}.(2). Let
$\maxLoc_l(R)$ be the set of maximal elements of the poset
$(\Loc_l(R), \ra )$. By the very definition of $\Loc_l(R)$ and by
 Lemma  \ref{a6Jul10}.(2), 
\begin{equation}\label{mADen1}
\maxLoc_l(R) = \{ S^{-1}R \, | \, S\in \maxDen_l(R) \}= \{ Q_l(R/
\ga ) \, | \, \ga \in \assmaxDen_l(R)\}.
\end{equation}

{\it Definition}. Each element  of $\maxLoc_l(R)$ is called a {\em
maximal left localization ring}   (or a {\em maximal left quotient
ring}) of the ring $R$.

\begin{theorem}\label{15Nov10}
Let $S\in \maxDen_l(R)$, $A= S^{-1}R$, $A^*$ be the group of units
of the ring $A$; $\ga := \ass (S)$, $\pi_\ga :R\ra R/ \ga $, $
a\mapsto a+\ga$, and $\s_\ga : R\ra A$, $ r\mapsto \frac{r}{1}$.
Then
\begin{enumerate}
\item $S=S_\ga (R)$, $S= \pi_\ga^{-1} (S_0(R/\ga ))$, $ \pi_\ga
(S) = S_0(R/ \ga )$ and $A= S_0( R/\ga )^{-1} R/ \ga = Q_l(R/ \ga
)$. \item  $S_0(A) = A^*$ and $S_0(A) \cap (R/ \ga )= S_0( R/ \ga
)$. \item $S= \s_\ga^{-1}(A^*)$. \item $A^* = \langle \pi_\ga (S)
, \pi_\ga (S)^{-1} \rangle$, i.e. the group of units of the ring
$A$ is generated by the sets $\pi_\ga (S)$ and $\pi_\ga^{-1}(S):=
\{ \pi_\ga (s)^{-1} \, | \, s\in S\}$. \item $A^* = \{ \pi_\ga
(s)^{-1}\pi_\ga ( t) \, |\, s, t\in S\}$. \item $Q_l(A) = A$ and
$\Ass_l(A) = \{ 0\}$.     In particular, if $T\in \Den_l(A, 0)$
then  $T\subseteq A^*$.
\end{enumerate}
\end{theorem}

{\it Proof}. 1. Since $S=S_\ga (R)$ (by Theorem \ref{3Jul10}.(2)),
statement 1 follows from Lemma \ref{a6Jul10}.(2).

2. $S_0(A)\stackrel{{\rm st. 1}}{=}S_0(Q_l(R/ \ga ))\stackrel{{\rm
Thm}\, \ref{4Jul10}.(1)}{=}Q_l(R/\ga )^* \stackrel{{\rm st.
1}}{=}A^*$ and $S_0(A)\cap ( R/\ga )\stackrel{{\rm st.
1}}{=}S_0(Q_l(R/ \ga ))\cap (R/\ga )\stackrel{{\rm Thm}
\ref{4Jul10}.(1)}{=}S_0(R/\ga )$.

3. $S= S_\ga (R)  \stackrel{{\rm Pr} \, \ref{8Jul10}.(1)}{=}
\s_\ga^{-1} (S_0 (S_\ga (R)^{-1}R))= \s_\ga^{-1} (S_0 (Q_l(R/\ga
))) \stackrel{{\rm Thm} \, \ref{4Jul10}.(1)}{=}\s_\ga^{-1}
(Q_l(R/\ga )^*) \stackrel{{\rm st. 1}}{=}\s_\ga^{-1} (A^*)$.

4. $A^*\stackrel{{\rm st. 1}}{=}Q_l(R/ \ga )^*\stackrel{{\rm
Thm}\,  \ref{4Jul10}.(3)}{=}\langle S_0(R/ \ga ), S_0(R/\ga
)^{-1}\rangle
 \stackrel{{\rm st. 1}}{=} \langle \pi_\ga (S) , \pi_\ga
 (S)^{-1}\rangle$.

5. $A^*\stackrel{{\rm st. 1}}{=}Q_l(R/ \ga )^*\stackrel{{\rm
Thm}\,  \ref{4Jul10}.(3)}{=}\{ p^{-1}q\, | \, p,q\in S_0(R/\ga )
\} = \{ \pi_\ga (s)^{-1} \pi_\ga (t)\, | \, s,t\in S\} $ since
$\pi_\ga (S) = S_0(R/\ga )$, by statement 1.

6. $Q_l(A)\stackrel{{\rm st. 1}}{=}Q_l(Q_l(R/ \ga ))\stackrel{{\rm
Thm}\,  \ref{4Jul10}.(4)}{=}Q_l(R/ \ga )\stackrel{{\rm st. 1}}{=}
A$.  Since $S_0(A)\stackrel{{\rm st. 1}}{=}S_0(Q_l(R/ \ga
))\stackrel{{\rm Thm}\, \ref{4Jul10}.(1)}{=}Q_l(R/ \ga )^* = A^*$,
$T\subseteq A^*$. The fact that $\Ass_l(A) = \{ 0\}$ follows from
Proposition \ref{A21Nov10}.(2) and the maximality of $A$. $\Box $

$\noindent $

The next theorem is a criterion of when a ring $A\in \Loc_l(R, \ga
)$ is equal to $Q_\ga (R)$.

\begin{theorem}\label{25Nov10}
Let $A\in \Loc_l(R, \ga )$, i.e. $A=S^{-1}R$ for some $S\in
\Den_l( R, \ga )$. Then $A= Q_\ga (R)$ iff $Q_l(A) = A$.
\end{theorem}

{\it Proof}. If $A= Q_\ga (R)$ then $A= Q_l( R / \ga )$, by Lemma
\ref{a6Jul10}.(2), and so $Q_l( A) = Q_l( Q_l( R/ \ga )) = Q_l( R/
\ga ) = A$, by Theorem \ref{4Jul10}.(4).

If $A\neq Q_\ga (R)$ then   the monomorphism $A\ra Q_\ga (R)$,
$a\mapsto a$, is not surjective. Let $S_1:= S$ and $S_2:= S_\ga
(R)$. Then $S_1\subseteq S_2$ and, by Lemma \ref{a5Jul10}.(2),
$Q_\ga (R) \simeq \pi_1(S_2)^{-1} (S_1^{-1}R) \simeq
\pi_1(S_2)^{-1}A$ where $\pi_1: R\ra R/ \ga$, $a\mapsto a+\ga$.
Therefore, $A\neq Q_l(A)$. $\Box $

 $\noindent $

{\bf Left localization maximal rings}. We introduce a new class of
rings, the left localization maximal rings, which turn out to be
precisely the class of maximal left quotient rings of all rings.
As a result, we have a characterization of the maximal left
quotient rings of a ring (Theorem  \ref{21Nov10}).

 $\noindent $

 {\it Definition}. A ring $A$ is
called a {\em left localization maximal ring} if $A= Q_l(A)$ and
$\Ass_l(A) = \{ 0\}$. A ring $A$ is called a {\em right
localization maximal ring} if $A= Q_r(A)$ and $\Ass_r(A) = \{
0\}$. A ring $A$ which is a left and right localization maximal
ring is called a {\em left and right localization maximal ring}
(i.e. $Q_l(A) =A=Q_r(A)$ and $\Ass_l(A) =\Ass_r(A) = \{ 0\}$).

$\noindent $

{\it Example}. Let $A$ be a simple ring. Then $Q_l(A)$ is a left
localization maximal  ring and $Q_r(A)$ is a right localization
maximal ring.

$\noindent $

{\it Example}. A division ring is a (left and right) localization
maximal ring. More generally, a simple artinian algebra (i.e. the
matrix algebra over a division ring) is a (left and right)
localization maximal ring.

$\noindent $

The next theorem is a criterion of  when a left quotient ring of a
ring is a maximal left quotient ring of the ring.

\begin{theorem}\label{21Nov10}
Let  a ring $A$ be a left localization of a ring $R$, i.e. $A\in
\Loc_l(R, \ga )$ for some $\ga \in \Ass_l( R)$. Then $A\in
\maxLoc_l( R)$ iff $Q_l( A) = A$ and  $\Ass_l(A) = \{ 0\}$, i.e.
$A$ is a left localization maximal ring.
\end{theorem}

{\it Proof}. $(\Rightarrow )$ Theorem \ref{15Nov10}.(6).

$(\Leftarrow )$ Proposition \ref{A21Nov10}. $\Box $

$\noindent $

The next corollary is a criterion of when $S_{\ga_1}(R)\subseteq
S_{\ga_2}(R)$ where $\ga_1, \ga_2\in \Ass_l(R)$.

\begin{corollary}\label{a25Nov10}
Let $\ga_1, \ga_2\in \Ass_l(R)$ and $\s_i : R\ra Q_{\ga_i}(R)$,
$r\mapsto r/1$, for $i=1,2$. Then $S_{\ga_1}(R) \subseteq
S_{\ga_2}(R)$ iff $\ga_1\subseteq \ga_2$ and $\s_2(S_{\ga_1}(R))
\subseteq  Q_{\ga_2}(R)^*$.
\end{corollary}

{\it Proof}. $(\Rightarrow )$ By Lemma \ref{a5Jul10}.(2a),
$\ga_1\subseteq \ga_2$ and, by Lemma \ref{a5Jul10}.(2b),
$\s_2(S_{\ga_1}(R))\subseteq     Q_{\ga_2}(R)^*$.

$(\Leftarrow )$ Let $S_i:= S_{\ga_i}(R)$ and $Q_i:= Q_{\ga_i}(R)$
for $i=1,2$. Let $Q'$ be the subring  of $Q_2$ generated by $R/
\ga_2$ and $\s_2(S_1)^{\pm 1}$. Since $S_1\in \Den_l(R, \ga_1)$, $
\ga_1\subseteq \ga_2$ and $\s_2(S_1) \subseteq Q_2^*$, every
element of $Q'$ has the form $\s_2(s_1)^{-1} \s_2(r)$ for some
elements $s_1\in S_1$ and $ r\in R$. By the universal property of
$Q_1= S_1^{-1}R$,  there exists a ring homomorphism $Q_1\ra Q_2$
and so we have the commutative diagram of ring homomorphisms:
$$
\xymatrix{ R\ar[r]\ar[d]^{=}  & R/ \ga_1 \ar[r]\ar[d] & Q_1 \ar[d]  \\
R\ar[r]  & R/ \ga_2 \ar[r] & Q_2 \, . }
$$
Since $S_i= \s^{-1}_i (Q_i^*)$ for $i=1,2$ (Lemma
\ref{a6Jul10}.(2) and Theorem \ref{4Jul10}), using the commutative
diagram we have the inclusion $S_1\subseteq S_2$. $\Box $


\begin{theorem}\label{C15Nov10}
Let $S\in \maxDen_l(R)$, $A=S^{-1}R$, $\ga = \ass (S)$ and $
\s_\ga : R\ra A$, $ r\mapsto \frac{r}{1}$. Then the following
statements are equivalent.
\begin{enumerate}
\item $A$ is a semi-simple ring. \item $Q_{cl} (R/ \ga ) $ exists
and is a semi-simple ring.
\end{enumerate}
If one of these conditions holds then $A=Q_{cl}(R)$ and $S=
\s_\ga^{-1} (Q_{cl}(R)^*)$.
\end{theorem}

{\it Proof}. 1. Since $A= Q_l(R/ \ga )$ and $S= \s_\ga^{-1}(A^*)$,
by Theorem \ref{15Nov10}.(1,3), the results follow from Theorem
\ref{5Jul10}. $\Box $

$\noindent $

{\bf The ideal $\gp (S)$ and the set $\CJ_l(R, S)$}.
\begin{proposition}\label{c14Nov10}
Let $S_1, S_2\in \Ore_l(R)$ be such that $\ga_1:= \ass (S_1)
\subseteq \ga_2:= \ass (S_2)$. Then
\begin{enumerate}
\item $S_1S_2\in \Ore_l(R, \ga )$  such that $ \ga_2\subseteq \ga
:= \ass (S_1S_2)$ where $S_1S_2:= \langle S_1, S_2 \rangle$ is the
sub-semigroup of $(R, \cdot )$ generated by $S_1$ and $S_2$. \item
If, in addition, $S_1, S_2, S\in \Den_l(R)$ then for each $i=1,2$
there is the ring homomorphism $S_i^{-1}R\ra S^{-1}R$,
$s^{-1}r\mapsto s^{-1}r$, with kernel $S_i^{-1}\ga$.
\end{enumerate}
\end{proposition}

{\it Proof}. 1. {\em Step 1}: $0\not\in S:=\langle S_1,
S_2\rangle$. Suppose that $0\in S$, we seek a contradiction. Then
$s_1s_2\cdots s_n =0$ for some elements $s_i\in S_1\cup S_2$ and
$n\geq 1$. Clearly, $n\geq 2$. We may assume that $n$ is the least
possible. Then, by the minimality of $n$, either all $s_{even}\in
S_1$ and $s_{odd}\in S_2$ or otherwise $s_{even}\in S_2$ and
$s_{odd}\in S_1$. If $s_1\in S_1$ then $s_2\in S_2$ and
$s_2s_3\cdots s_n\in \ga_1\subseteq \ga_2$, hence
$(s_2's_2)s_2\cdots s_n =0$ for some element $s_2'\in S_2$. This
contradicts to the minimality of $n$ since $s_2's_2\in S_2$. So,
$s_1\in S_2$. If $s_n \in S_1$ then $s_1's_1s_2\cdots s_n =0$ for
some element $s_1'\in S_1$. This is not possible by the previous
case. Therefore $s_1, s_n\in S_2$. Then $ (s_n's_1)s_2\cdots
s_{n-1} =0$ for some element $s_n'\in S_2$. This contradicts to
the minimality of $n$ since $s_n's_1\in S_2$. Therefore, $0\not\in
S$.

 {\em Step 2}: $S\in \Ore_l(R)$. We have to show that for any
 elements $s\in S$ and $r\in R$ there exist elements $ s'\in S$
 and $r'\in R$ such that $s'r=r's$. To prove this we use induction
 on $n$ where $s= s_1s_2\cdots s_n$, $s_i\in S_1\cup S_2$. The
 result is obvious when $n=1$ since $ S_1,S_2\in \Ore_l(R)$.
 Suppose that $n>1$ and the result is true for all $n'<n$. Fix
 $t_1\in S_1\cup S_2$ such that $t_1r=r_1s_n$ for some element
 $r_1\in R$. By induction, $t_2r_1=r's_1\cdots s_{n-1}$ for some $t_2\in
 S$ and $r'\in R$. Then $t_2t_1r= t_2 r_1s_n = r' s_1\cdots s_n$.
  It suffices to take $s'=t_2t_1$. By Step 2, $\ga $ is an ideal of
 the ring $R$ such that $\ga \neq R$, by Step 1. Clearly, $\ga_2\subseteq
 \ga$.

2. Statement 2 follows from statement 1 and Lemma
\ref{a5Jul10}.(2a). $\Box $

$\noindent $

In order to give an answer to the question of when  an Ore set $S$
 of a ring is a denominator set  or, more generally, when the image of $S$
 is a denominator set in a factor ring,    we introduce the ideal $\gp (S)$.

 Let $\mW$ be the family  of all {\em ordinal} numbers. For each left
 Ore set $S\in \Ore_l(R)$, we attach the ideal of the ring $R$,
\begin{equation}\label{paup}
\gp (S) :=\gp (S)_l:= \bigcup_{\alpha \in \mW}\gp_\alpha
\end{equation}
where the ideals $\gp_\alpha$ of the ring $R$ are defined
recursively as follows: $\gp_1 :=\ass_l (S, R) +\ass_r(S, R)R$
where $\ass_l(S, R):=\{ r\in R\, | \, sr=0$ for some $s=s(r)\in
S\}$ and $\ass_r(S, R):= \{ r\in R\, | \, rs=0$ for some $s=s(r)
\in S\}$. Note that $\ass_l(S, R)$ is an ideal of the ring $R$
since $S\in \Ore_l(R)$ and  $\ass_r(S, R)$ is a right ideal of the
ring $R$. 
\begin{equation}\label{paup1}
\gp_{\alpha +1}:= \pi_\alpha^{-1}(\ass_l(\pi_\alpha (S), R/
\gp_\alpha ) +\ass_r(\pi_\alpha (S), R/\gp_\alpha )\cdot R/
\gp_\alpha)
\end{equation}
 where $\pi_\alpha : R\ra R/\gp_\alpha$, $a\mapsto
a+\gp_\alpha$. If $\alpha$ is a limit ordinal then
\begin{equation}\label{paup2}
\gp_\alpha = \bigcup_{\beta <\alpha} \gp_\beta .
\end{equation}
By the very definition, if $\alpha \leq \beta$ then $ \gp_\alpha
\leq \gp_\beta$. In a similar fashion, for a {\em right} Ore set
$S\in \Ore_r(R)$, we can define the ideal $\gp (S)_r$. For each
left Ore set $S\in \Ore_l(R)$, let $$\CJ_l(R, S):=\{ \ga\, | \,
\ga\; {\rm is \; an \; ideal\; of\;} R\; {\rm  such \; that } \;
\pi_\ga (S)\in \Den_l(R/ \ga , 0)\}$$ where $\pi_\alpha : R\ra R/
\ga$, $a\mapsto a+\ga$.

\begin{theorem}\label{A15Nov10}
Let $S\in \Ore_l(R)$ and $\ga = \ass (S)$. Then
\begin{enumerate}
\item $\CJ_l(R,S) \neq \emptyset$ iff $\gp (S) \neq R$. \item
 If $\CJ_l(R,S) \neq \emptyset$ then $\gp (S)$ is
 the least (with respect to inclusion) element of $\CJ_l(R, S)$.
 \item If, in addition, $S\in \Den_l(R, \ga )$ then $\gp (S) =
 \ga$.
 \item If the ring $R$ satisfies the ascending chain condition on
 annihilator right ideals then $S\in \Den_l(R, \ga )$ and $\gp (S)
 = \ga$.
\end{enumerate}
\end{theorem}

{\it Proof}. 1. $(\Rightarrow )$ Suppose that $\CJ_l(R,S) \neq
\emptyset$. Fix $\ga \in \CJ_l(R,S)$.  Clearly, $\ass_l(S, R),
\ass_r(S,R)\subseteq \ga$, hence $\gp_1\subseteq \ga$. If
$\gp_\alpha \subseteq \ga$ then $\gp_{\alpha +1} \subseteq \ga$
since $sa, bt\in \gp_\alpha \subseteq \ga$ for some $s,t\in S$ and
$ a,b\in R$ implies $a,b\in \ga$. If $\alpha$ is a limit ordinal
and $\gp_\beta \subseteq \ga$ for all $\beta <\alpha$ then
$\gp_\alpha = \bigcup_{\beta < \alpha} \gp_\beta \subseteq \ga$.
Therefore, $\gp (S)\subseteq \ga$, and so $\gp (S) \neq R$.

$(\Leftarrow )$ Suppose that $\gp:= \gp (S) \neq R$. We claim that
$\pi_\gp (S) \in \Den_l(R/\gp , 0)$ where $\pi_\gp : R\ra R/\gp$,
$a\mapsto a+\ga$. If $\pi_\gp (s) \pi_\gp (a) =0$ and $\pi_\gp (b)
\pi_\gp (t)=0$ in $R/\gp$ for some elements $s,t\in S$ and $a,b\in
R$ then $sa, bt\in \gp =\bigcup_{\alpha \in \mW}\gp_\alpha$, i.e.
$sa, bt\in \gp_\alpha$ for some ordinal number $\alpha$.
Therefore, $a,b\in \gp_{\alpha +1} \subseteq \gp$, i.e. $\pi_\gp
(a) =\pi_\gp (b) =0$. This proves that $\ass_l( \pi_\gp (S), R/\gp
)=0$ and $\ass_r( \pi_\gp (S), R/ \gp )=0$. To finish the proof of
the fact that $\pi_\gp (S)\in \Den_l(R/\gp , 0)$ we have to show
that $\pi_\gp (S) \in \Ore_l(R/\gp )$, but his is obvious as
epimorphisms respect left Ore sets provided their images are
multiplicatively closed sets which is the case for $\pi_\gp (S)$
as $\gp (S) \neq R$.

2. In statement 1 we proved that $\gp (S)\subseteq \ga$ for all
ideals $\ga \in \CJ_l(R,S)$ and $\gp (S) \in \CJ_l(R, S)$, i.e.
$\gp (S)$ is the least element of the poset $(\CJ_l(R, S),
\subseteq )$.

3. Statement 3 is obvious.

4. By Lemma 1.1.2, \cite{Jategaonkar-book}, $S\in \Den_l(R, \ga )$
and so $\gp (S) = \ga$, by statement 3. $\Box $

$\noindent $

{\bf Localizable left Ore sets, a slight generalization of Ore's
method of localization}.

{\it Definition}. A left Ore set $S\in \Ore_l(R)$ is called {\em
localizable} if there exists a ring homomorphism $\v : R\ra R'$
where $R'$ is a ring such that, for all $s\in S$, $\v (s)$ is a
unit of the ring $R'$. The set of all localizable left Ore sets in
$R$ is denoted by $\ore_l(R)$. By Theorem \ref{A15Nov10},
\begin{equation}\label{orelSp}
\ore_l(R)=\{ S\in \Ore_l(R)\, | \, \gp (S) \neq R\}.
\end{equation}

\begin{corollary}\label{28Nov10}
Let $S$ be a localizable left Ore set in a ring $R$. Then there
exists an ordered pair $(Q, f)$ where $Q$ is a ring and $f: R\ra
Q$ is a ring homomorphism such that

(i)  for all $s\in S$, $f(s)$ is a unit in $Q$; 


 and if $(Q', f')$ is another
 pair satisfying the condition
 (i) then there is a unique ring homomorphism $h: Q\ra Q'$ such that
$f' = hf$. The ring $Q$ is unique up to isomorphism. The ring $Q$
is isomorphic to the left localization of the ring $R/ \gp (S)$ at
the left denominator set $\pi (S)\in \Den_l(R/ \gp (S), 0)$ where
the ideal $\gp (S)$ of the ring $R$ is defined in (\ref{paup2})
and $\pi : R\ra R/ \gp (S)$, $ a\mapsto a+\gp (S)$.
\end{corollary}

{\it Proof}. Recall that $\gp (S) = \bigcup_{\alpha \in\mW }
\gp_\alpha$. By induction on $\alpha$, all $\gp_\alpha \subseteq
\ker (f')$. Therefore, $\gp (S) \subseteq \ker (f')$. Notice that
$Q' \simeq \pi_{\ga}(S)^{-1}(R/ \ga )$ for some $\ga \in \CJ_l(R,
S)$. Let us identify this rings via the isomorphism above, then
$f'(R) = R/ \ga$. By Theorem \ref{A15Nov10}.(2), $\gp (S)\subseteq
\ga$, and so there is a natural ring homomorphism $\v : R/ \gp (S)
\ra \pi_{\ga} (S)^{-1}(R/ \ga )$ such that, for all $s\in S$, $\v
(s)$ is a unit. By the universal property of the ring
$\pi_{\gp(S)}(S)^{-1} (R/ \gp (S))$ there is a {\em unique}
homomorphism $h:\pi_{\gp (S)}(S)^{-1} (R/ \gp (S)) \ra Q'$ such
that $f'=hf$. This proves existence of the ring $Q$. It is unique
up to isomorphism as it follows  at once from the uniqueness of
$h$. $\Box $

$\noindent $

The set $\LDen_l( R) := \{ S_\ga (R) \, | \, \ga \in \Ass_l(R)\}$
is the set of largest left denominator sets of the ring $R$. Let
$\minLDen_l( R, \ga ) $ and $\minLoc_l(R, \ga )$ be the sets of
minimal elements of the posets $\LDen_l(R, \ga )$ and $\Loc_l( R,
\ga )$ respectively.

%



\section{The largest  quotient ring of a ring}\label{LQR}

In this section, the two-sided analogues of some results of
Sections \ref{LLQRR} and \ref{TMLQR} are proved. The proofs are
dropped in all cases where the two-sided analogues are direct
corollaries of the one-sided results. In particular, it is proved
existence of the largest  quotient ring of a ring  (Theorem
\ref{T3Jul10}) and that all Ore sets are localizable (Theorem
\ref{B15Nov10}). The results of this section are used in Section
\ref{EXEM}.

{\bf Notation}:

\begin{itemize}
\item $\Ore (R) := \Ore_l(R)\cap \Ore_r(R)= \{ S\, | \, S$ is a
left and right  Ore set in $R\}$; \item $\Den(R):=\Den_l(R)\cap
\Den_r(R)=\{ S\, | \, S$ is a left and right denominator set in
$R\}$. For each $S\in \Den (R)$, $\ass (S):=\ass_l (S)=\ass_r (S)$
is an ideal of the ring $R$  where $\ass_l(S):=\{ r\in R\, | \,
sr=0$ for some $s=s(r)\in S\}$ and $\ass_r(S):=\{ r\in R\, | \,
rs=0$ for some $s=s(r)\in S\}$; \item $\Loc (R):= \{
S^{-1}R=RS^{-1}\, | \, S\in \Den(R)\}$; \item $\Ass(R):= \{ \ass
(S)\, | \, S\in \Den(R)\}$; \item $\Den (R, \ga ) := \{ S\in \Den
(R)\, | \, \ass (S)=\ga \}$ where $\ga \in \Ass (R)$; \item
$S_\ga=S_\ga (R)$ is the {\em largest element} of the poset $(\Den
(R, \ga ), \subseteq )$ and $Q_\ga (R):=S_\ga^{-1} R=RS_\ga^{-1}$
is the {\em largest (two-sided)  quotient ring associated to}
$\ga$, $S_\ga $ exists (Theorem \ref{T3Jul10}.(2)); \item In
particular, $S_0=S_0(R)$ is the largest element of the poset
$(\Den (R, 0), \subseteq )$ and $Q (R):=S_0^{-1}R=RS_0^{-1}$ is
the {\em largest (two-sided) quotient ring} of $R$; \item $\Loc
(R, \ga ):= \{ S^{-1}R=RS^{-1}\, | \, S\in \Den (R, \ga )\}$.
\end{itemize}

{\it Remark}. Subscripts `l' and `r' indicate that  left and right
versions of a definition/concept are considered respectively.

$\noindent $

{\bf The largest quotient ring of a ring}.

\begin{theorem}\label{T3Jul10}
\begin{enumerate}
\item For each $\ga \in \Ass(R)$, the set $\Den(R, \ga )$   is an
ordered abelian semigroup $(S_1S_2= S_2S_1$, and $S_1\subseteq
S_2$ implies $S_1S_3\subseteq S_2S_3$) where the product $S_1S_2=
\langle S_1, S_2\rangle$ is the multiplicative subsemigroup of
$(R, \cdot )$ generated by $S_1$ and $S_2$.  \item  $S_\ga :=
S_\ga (R) := \bigcup_{S\in \Den(R, \ga )}S$ is the largest element
(w.r.t. $\subseteq $) in $\Den(R, \ga )$. The set $S_\ga$ is
called the {\em largest  denominator set} associated to $\ga$.
\item Let $S_i\in \Den(R,\ga )$, $i\in I$, where $I$ is an
arbitrary non-empty set. Then 
\begin{equation}\label{SiIgen}
\langle S_i\, | \, i\in I\rangle:=\bigcup_{\emptyset \neq
J\subseteq I, |J|<\infty}\prod_{j\in J}S_j\in \Den(R, \ga )
\end{equation}
the  denominator set generated by the denominators sets $S_i$, it
is the least upper bound of the set $\{ S_i\}_{i\in I}$ in
$\Den(R, \ga )$, i.e. $\langle S_i\, | \, i\in I\rangle=
\bigvee_{i\in I}S_i$.
\end{enumerate}
\end{theorem}

In Section \ref{EXEM}, we will see that for the algebra $\mI_1$ of
polynomial integro-differential operators the set $S_0(\mI_1)$
(resp. the ring $Q(\mI_1)$) is tiny comparing to the sets
$S_{l,0}(\mI_1)$ and $S_{r, 0}(\mI_1)$ (resp. to the rings
$Q_l(\mI_1)$ and $Q_r(\mI_1)$).

\begin{corollary}\label{Ta7Jul10}
The abelian monoid $\Den(R, 0)$ is a complete lattice such that
 $S_1S_2=S_1\bigvee S_1$ and $\bigwedge_{i\in I}S_i=\bigvee_{S\in
 \Den(R, 0), S\subseteq \cap_{i\in I}S_i}S_i$ where all $S_i\in
 \Den(R, 0)$.
\end{corollary}

Clearly, $\bigwedge_{i\in I}S_i$ is the largest element of the set
$\{ S\, | \, S\in \Den(R,0), S\subseteq \bigcap_{i\in I}S_i\}$.

\begin{corollary}\label{T7Jul10}
\begin{enumerate}
\item Let $R$ be a ring. Each ring automorphism $\s\in \Aut(R)$ of
the ring $R$ has the unique extension $\s\in \Aut(Q(R))$ to an
automorphism of the ring $Q(R)$ given by the rule $\s (s^{-1} r) =
\s (s)^{-1}\s (r)$ where $s\in S_0(R)$ and $r\in R$. \item The
group $\Aut(R)$ is a subgroup of the group $\Aut(Q(R))$. Moreover,
$\Aut (R) = \{ \tau \in \Aut(Q(R))\, | \, \tau (R) = R\}$.
\end{enumerate}
\end{corollary}

\begin{theorem}\label{T4Jul10}
\begin{enumerate}
\item $ S_0 (Q(R))= Q(R)^*$ {\em and} $S_0(Q(R))\cap R= S_0(R)$.
 \item $Q(R)^*= \langle S_0(R), S_0(R)^{-1}\rangle$, {\em i.e. the
 group of units of the ring $Q(R)$ is generated by the sets
 $S_0(R)$ and} $S_0(R)^{-1}:= \{ s^{-1} \, | \, s\in S_0(R)\}$.
 \item $Q(R)^* = \{ s^{-1}t\, | \, s,t\in S_0(R)\}= \{ ts^{-1}\, | \, s,t\in S_0(R)\}$.
 \item $Q(Q(R))=Q(R)$.
\end{enumerate}
\end{theorem}

{\it Proof}. 1--3. It is obvious that $G:=\langle S_0(R),
S_0(R)^{-1}\rangle\subseteq Q(R)^*\subseteq S_0(Q(R))$. Applying
Proposition \ref{A4Jul10}.(2) and its right analogue  in the
situation where $S=S_0(R)$ and $T = S_0(Q(R))$ we see that
$$S_0(R)\subseteq T':= S_0(Q(R))\cap R\in \Den_l(R, 0)\cap \Den_r(R, 0),$$
and so $T'\in \Den (R, 0)$ and $S_0(R) = T'$, by the maximality of
$S_0(R)$. Let $q\in S_0(Q(R))$. Then $q= s^{-1} t=t_1s_1^{-1}$ for
some elements $s, s_1\in S_0(R)$, $t,t_1\in R$,  $t = sq\in
S_0(Q(R))\cap R= S_0(R)$ and $t_1 = qs_1\in S_0(Q(R))\cap R=
S_0(R)$. Therefore, $S_0(Q(R)) \subseteq \{ s^{-1} t\, | \, s,t\in
S_0(R)\} \subseteq G$ and $S_0(Q(R)) \subseteq \{ t s^{-1}\, | \,
s,t\in S_0(R)\} \subseteq G$, and so $G= Q(R)^* = S_0(Q(R)) = \{
s^{-1} t\, | \, s,t \in S_0(R)\} = \{ t s^{-1}\, | \, s,t \in
S_0(R)\}$. This proves statements 1--3.

4. Statement 4 follows from statement 1. $\Box $


\begin{theorem}\label{T5Jul10}
Let $R$ be a ring and $\CC_R$ be the set of regular elements of
the ring $R$ (i.e. the set of non-zero-divisors).  Then the
following statements are equivalent.
\begin{enumerate}
\item $Q(R)$ is a semi-simple ring. \item The rings $Q_{l, cl}(R)$
and $Q_{r, cl}(R)$ exist and are semi-simple. \item The rings
$Q_{l}(R)$ and $Q_{r}(R)$  are semi-simple.
\end{enumerate}
If one of the equivalent conditions holds then $S_0(R)= \CC_R$ and
$Q(R) \simeq Q_{l, cl}(R)\simeq Q_{r, cl}(R)\simeq Q_{l}(R)\simeq
Q_{r}(R)$.
\end{theorem}

{\it Proof}. $(1\Rightarrow 3)$ Since $Q(R)\subseteq Q_l(R)$,
$Q(R)\subseteq Q_r(R)$ and the ring $Q(R)$ is a semi-simple, by
Lemma \ref{a5Jul10}.(2),  the inclusions are, in fact, equalities
(since the regular elements  of any semi-simple ring are the
units). Therefore, the rings $Q_{l}(R)$ and $Q_{r}(R)$  are
semi-simple.

$(2\Leftrightarrow 3)$ Statements 2 and 3 are equivalent, by
Theorem  \ref{5Jul10}. Moreover, $Q_{l, cl}(R)\simeq Q_{r,
cl}(R)\simeq Q_{l}(R)\simeq Q_{r}(R)$.

$(2\Rightarrow 1)$ Recall that $Q_{l, cl}(R)=\CC_R^{-1}R$ and
$Q_{r, cl}(R)=R\CC_R^{-1}$ where $\CC_R$ is the set of regular
elements of the ring $R$. Since the rings  $Q_{l, cl}(R)$ and
$Q_{r, cl}(R)$  exist, they are $R$-isomorphic. Therefore, $S_0(R)
= \CC_R$ and $Q(R)\simeq Q_{l, cl}(R)\simeq Q_{r, cl}(R)$. In
particular, $Q(R)$
 is a semi-simple ring. $\Box $


\begin{proposition}\label{TA5Jul10}
If the ring $Q_l(R)$ is a left artinian ring  and the  ring
$Q_r(R)$ is a right artinian ring then $S_0(R)= \CC_R$ and $Q(R)
\simeq Q_{l, cl}(R)\simeq Q_{r, cl}(R)\simeq Q_{l}(R)\simeq
Q_{r}(R)$.
\end{proposition}

{\it Proof}. By Corollary \ref{A5Jul10}, $S_{l,0}(R) = \CC_R$ and
$Q_{l, cl}(R)\simeq Q_l(R)$. By the right version of Corollary
\ref{A5Jul10}, $S_{r,0}(R) = \CC_R$ and $Q_{r, cl}(R)\simeq
Q_r(R)$. Then $S_0(R) = \CC_R$ and $Q(R) \simeq Q_{l, cl}(R)\simeq
Q_{r, cl}(R)\simeq Q_{l}(R)\simeq Q_{r}(R)$.  $\Box $

$\noindent $

{\bf The maximal quotients rings of a ring}.

\begin{lemma}\label{Ta6Jul10}
Let $S\in \Den(R, \ga )$, $\pi : R\ra R/\ga$, $a\mapsto a+\ga$,
and  $\s : R\ra S^{-1}R$, $ r\mapsto r/1$.
\begin{enumerate}
\item  Let $T\in \Den(S^{-1}R, 0)$ be such that $\pi (S), \pi
(S^{-1})\subseteq T$. Then $T':= \s^{-1} (T) \in \Den(R, \ga )$,
$T'$ is $S$-saturated, $T=\{ s^{-1}t'\, | \, s\in S, t'\in T'\}$,
and $S^{-1} R\subseteq T'^{-1}R= T^{-1}R$. \item $\pi^{-1} (S_0(R/
\ga )) = S_\ga (R)$, $ \pi (S_\ga (R)) = S_0(R/\ga ))$ and $ Q_\ga
(R)= Q(R/\ga )$.
\end{enumerate}
\end{lemma}

\begin{proposition}\label{T8Jul10}
Let $S\in \Den(R, \ga )$; $\pi :R\ra R/ \ga$, $a\ra \oa = a+\ga$;
 $\s : R\ra S^{-1}R, r\ra r/1$; $G:=\langle \pi (S), \pi
 (S)^{-1}\rangle \subseteq (S^{-1}R)^*$ (i.e. $G$ is the subgroup
 of the group $(S^{-1}R)^*$ of units of the ring $S^{-1}R$
 generated by $\pi (S)^{\pm 1}$). Let $[\s^{-1} (G), S_\ga (R)]:=
 \{ S\in \Den (R, \ga )\, | \,  \s^{-1} (G)\subseteq S \subseteq
 S_\ga (R)\}$.
\begin{enumerate}
\item Let  $ [\s^{-1} (G), S_\ga (R)]_{S-{\rm
 com}}:= \{ S_1\in [\s^{-1} (G), S_\ga (R)]\, | \, \s^{-1} (G\pi (S_1))=S_1\}$ and
 $[ G, S_0(S^{-1}R)]:= \{ T \in \Den(S^{-1}R, 0)\, | \, G\subseteq T\subseteq S_0(S^{-1}R)\}$. Then the map
$$[\s^{-1}(G), S_\ga (R)]_{S-{\rm com}}\ra [ G, S_0(S^{-1}R)], \;\;
S_1\mapsto \ttS_1:= G\pi (S_1), $$ is an isomorphism of posets and
abelian monoids with the inverse map $T\mapsto \s^{-1}(T)$ where
$G\pi (S_1)$ is the multiplicative monoid generated by $G$ and
$\pi (S_1)$ in $S^{-1}R$. In particular,
$$G\pi (S_\ga (R)) = S_0(S^{-1}R), \;\; S_\ga(R) = \s^{-1}
(S_0(S^{-1}R)), \;\; S_\ga (R)^{-1}R \simeq S_\ga (R)^{-1}R\simeq
Q(R/\ga ),$$ the monoid $[\s^{-1}(G), S_\ga (R)]_{S-{\rm
 com}}$ is a complete lattice  (since $[ G, S_0(S^{-1}R)]$ is a complete  lattice, as an interval
 of the complete lattice $\Den(S^{-1}R, 0)$, Corollary \ref{Ta7Jul10}), and the
 map $S_1\mapsto \widetilde{S}_1$ is a lattice isomorphism.
 \item Consider the interval $[G\cap (R/ \ga ) , S_0(R/ \ga )]$ in
 $\Den(R/ \ga  , 0)$. Let $[G\cap (R/ \ga ) , S_0(R/ \ga
 )]_{G-{\rm com}}:=\{ T\in [G\cap (R/ \ga ) , S_0(R/ \ga )]\, | \,
 GT\cap (R/ \ga ) =T\}$. Then $[G\cap (R/ \ga ) , S_0(R/\ga )]_{G-{\rm com}} \subseteq
 \Den(S^{-1}R, 0)$ and  the map
 $$[G\cap (R/ \ga ) , S_0(R/ \ga
 )]_{G-{\rm com}}\ra [G, S_0(S^{-1}R)], \;\; T\mapsto GT,$$
 is an isomorphism of posets and abelian monoids with the inverse
 map $T'\mapsto T'\cap (R/ \ga )$ where $GT$ is the product in $\Den(S^{-1}R, 0)$. In particular,
 the monoid $[G\cap (R/ \ga ) , S_0(R/ \ga
 )]_{G-{\rm com}}$ is a complete lattice.
  \item The map $$[ \s^{-1} (G), S_\ga (R)]_{S-{\rm com}}\ra [G\cap (R/ \ga ) , S_0(R/ \ga
 )]_{G-{\rm com}}, \;\; S_1\mapsto G\pi (S_1)\cap (R/ \ga ),$$
is an isomorphism of posets and abelian monoids with the inverse
 map $S'\mapsto\pi^{-1}(S')$.
\end{enumerate}
\end{proposition}

\begin{lemma}\label{Tb14Nov10}
\begin{enumerate}
\item  For each element $S\in \Den(R)$, there exists a maximal
element in the poset $([S, \cdot ), \subseteq )$. \item The set
$(\maxDen(R), \subseteq )$ of maximal elements of the poset
$(\Den(R), \subseteq )$ is a non-empty set.
\end{enumerate}
\end{lemma}


{\it Definition}. An element $S$ of the set $\maxDen (R)$ is
called a {\em maximal  denominator set} of $R$ and the ring
$S^{-1}R=RS^{-1}$ is called a {\em maximal  quotient ring} of $R$.

\begin{proposition}\label{TA21Nov10}
Let $\ga \in \Ass(R)$, $Q:= Q_\ga (R)$, $Q^*$ be the group of
units of the ring $Q$ and  $\s : R\ra Q_\ga (R)$, $ r\mapsto
\frac{r}{1}$. Let $T\in \Den(Q, \gb )$ where $\gb \in \Ass(Q)$ and
$Q^*T$ be the multiplicative sub-semigroup of $(Q, \cdot )$
generated by $Q^*$ and $T$. Then
\begin{enumerate}
\item $Q^*T\in \Den(Q, \gb )$. \item If, in addition, $Q^*
\subseteq T$ (eg, $Q^*T$ from statement 1) then
\begin{enumerate}
\item $T':= \s^{-1} (T) \in \Den(R, \gb')$ where $\gb':= \s^{-1}
(\gb )\supseteq \ga$, $S_\ga (R) \subseteq T'$, $T= Q^* \s (T')$
(i.e. the monoid $T$ is generated by $Q^*$ and $\s (T')$) and
$T'^{-1}R = T^{-1}Q$ (i.e. the natural ring monomorphism
$T'^{-1}R\ra T^{-1}Q$, $ t^{-1}r\mapsto t^{-1}r$, is an
isomorphism). \item $S_\ga (R) \subseteq S_{\gb'}(R)$ and
$S_{\gb'} (R) = \s^{-1} (S_\gb (Q))$. \item $Q_{\gb'}(R) = Q(
Q_\ga (R)/ \gb )$, i.e. the natural ring monomorphism
$Q_{\gb'}(R)\ra Q( Q_\ga (R)/ \gb )$, $ t^{-1} r \mapsto t^{-1}r$,
is an isomorphism.
\end{enumerate}
\end{enumerate}
\end{proposition}

Let $\maxAss (R)$ be the set of maximal elements of the poset
$(\Ass(R), \subseteq )$. It is a subset of the set

\begin{equation}\label{TmADen}
\assmaxDen (R) := \{ \ass (S) \, | \, S\in \maxDen (R) \}
\end{equation}
which  is a non-empty set, by Lemma  \ref{Tb14Nov10}.(2). Let
$\maxLoc (R)$ be the set of maximal elements of the poset $(\Loc
(R), \ra )$. By the very definition of $\Loc (R)$ and by
 Lemma  \ref{Ta6Jul10}.(2), 
\begin{equation}\label{TmADen1}
\maxLoc (R) = \{ S^{-1}R \, | \, S\in \maxDen (R) \}= \{ Q (R/ \ga
) \, | \, \ga \in \assmaxDen (R)\}.
\end{equation}

{\it Definition}. Each element  of $\maxLoc (R)$ is called a {\em
maximal  localization ring}   (or a {\em maximal   quotient ring})
of the ring $R$.

\begin{theorem}\label{T15Nov10}
Let $S\in \maxDen (R)$, $A= S^{-1}R$, $A^*$ be the group of units
of the ring $A$; $\ga := \ass (S)$, $\pi_\ga :R\ra R/ \ga $, $
a\mapsto a+\ga$, and $\s_\ga : R\ra A$, $ r\mapsto \frac{r}{1}$.
Then
\begin{enumerate}
\item $S=S_\ga (R)$, $S= \pi_\ga^{-1} (S_0(R/\ga ))$, $ \pi_\ga
(S) = S_0(R/ \ga )$ and $A= S_0( R/\ga )^{-1} (R/ \ga )=(R/\ga
)S_0( R/\ga )^{-1} = Q (R/ \ga )$. \item  $S_0(A) = A^*$ and
$S_0(A) \cap (R/ \ga )= S_0( R/ \ga )$. \item $S=
\s_\ga^{-1}(A^*)$. \item $A^* = \langle \pi_\ga (S) , \pi_\ga
(S)^{-1} \rangle$, i.e. the group of units of the ring $A$ is
generated by the sets $\pi_\ga (S)$ and $\pi_\ga^{-1}(S):= \{
\pi_\ga (s)^{-1} \, | \, s\in S\}$. \item $A^* = \{ \pi_\ga
(s)^{-1}\pi_\ga ( t) \, |\, s, t\in S\}$. \item $Q (A) = A$ and
$\Ass(A) = \{ 0\}$.     In particular, if $T\in \Den (A, 0)$ then
$T\subseteq A^*$.
\end{enumerate}
\end{theorem}


The next theorem is a criterion of when a ring $A\in \Loc (R, \ga
)$ is equal to $Q_\ga (R)$.

\begin{theorem}\label{T25Nov10}
Let $A\in \Loc (R, \ga )$, i.e. $A=S^{-1}R$ for some $S\in \Den(
R, \ga )$. Then $A= Q_\ga (R)$ iff $Q (A) = A$.
\end{theorem}

{\bf Localization maximal rings}. We introduce a new class of
rings, the  localization maximal rings, which turn out to be
precisely the class of maximal  quotient rings of all rings. As a
result, we have a characterization of the maximal  quotient rings
of a ring (Theorem  \ref{T21Nov10}).

 $\noindent $

 {\it Definition}. A ring $A$ is
called a {\em localization maximal ring} if $A= Q (A)$ and
$\Ass(A) = \{ 0\}$. Clearly, a left and right localization maximal
ring  $A$  (i.e. $Q_l(A) =A=Q_r(A)$ and $\Ass(A) =\Ass_r(A) = \{
0\}$) is a localization maximal but vice versa, in general, as the
example of the algebra $\mI_1$ shows (see Section \ref{EXEM} for
details).


The next theorem is a criterion of  when a quotient ring of a ring
is a maximal  quotient ring of the ring.

\begin{theorem}\label{T21Nov10}
Let  a ring $A$ be a  localization of a ring $R$, i.e. $A\in \Loc
(R, \ga )$ for some $\ga \in \Ass( R)$. Then $A\in \maxLoc ( R)$
iff $Q ( A) = A$ and  $\Ass(A) = \{ 0\}$, i.e. $A$ is a
localization maximal ring.
\end{theorem}


The next corollary is a criterion of when $S_{\ga_1}(R)\subseteq
S_{\ga_2}(R)$ where $\ga_1, \ga_2\in \Ass(R)$.

\begin{corollary}\label{Ta25Nov10}
Let $\ga_1, \ga_2\in \Ass(R)$ and $\s_i : R\ra Q_{\ga_i}(R)$,
$r\mapsto r/1$, for $i=1,2$. Then $S_{\ga_1}(R) \subseteq
S_{\ga_2}(R)$ iff $\ga_1\subseteq \ga_2$ and $\s_2(S_{\ga_1}(R))
\subseteq  Q_{\ga_2}(R)^*$.
\end{corollary}

{\it Proof}. $(\Rightarrow )$ By Lemma \ref{a5Jul10}.(2a),
$\ga_1\subseteq \ga_2$ and, by Lemma \ref{a5Jul10}.(2b) and its
right analogue, $\s_2(S_{\ga_1}(R))\subseteq     Q_{\ga_2}(R)^*$.

$(\Leftarrow )$ Let $S_i:= S_{\ga_i}(R)$ and $Q_i:= Q_{\ga_i}(R)$
for $i=1,2$. Let $Q'$ be the subring  of $Q_2$ generated by $R/
\ga_2$ and $\s_2(S_1)^{\pm 1}$. Since $S_1\in \Den (R, \ga_1)$, $
\ga_1\subseteq \ga_2$ and $\s_2(S_1) \subseteq Q_2^*$, every
element of $Q'$ has the form $\s_2(s)^{-1}
\s_2(r)=\s_2(r')\s_2(s')^{-1}$ for some elements $s, s'\in S_1$
and $ r, r'\in R$. By the universal property of $Q_1=
S_1^{-1}R=RS_1^{-1}$, there exists a ring homomorphism $Q_1\ra
Q_2$ and so we have the commutative diagram of ring homomorphisms:
$$
\xymatrix{ R\ar[r]\ar[d]^{=}  & R/ \ga_1 \ar[r]\ar[d] & Q_1 \ar[d]  \\
R\ar[r]  & R/ \ga_2 \ar[r] & Q_2 \, . }
$$
Since $S_i= \s^{-1}_i (Q_i^*)$ for $i=1,2$ (Lemma
\ref{Ta6Jul10}.(2)), we have the inclusion $S_1\subseteq S_2$.
$\Box $

$\noindent $

{\bf All Ore sets are localizable, a slight generalization of
Ore's method of localization}.

{\it Definition}. An Ore set $S\in \Ore (R)$ is called {\em
localizable} if there exists a ring homomorphism $\v : R\ra R'$
where $R'$ is a ring such that, for all $s\in S$, $\v (s)$ is a
unit of the ring $R'$.

We will see that all Ore sets are localizable (Theorem
\ref{B15Nov10} and Corollary \ref{T28Nov10}).

For each right Ore set $S\in \Ore_r (R)$ of the ring $R$, let $\gp
(S)_r:= \bigcup_{\alpha \in \mW}\gp_\alpha'$ be the right version
of the ideal $\gp (S)_l$ defined in (\ref{paup}) where
\begin{equation}\label{rpaup1}
\gp_{\alpha +1}':= \pi_\alpha'^{-1} (R/ \gp_\alpha \cdot
\ass_l(\pi_\alpha' (S), R/ \gp_\alpha ) +\ass_r(\pi_\alpha' (S),
R/ \gp_\alpha )),
\end{equation}
$\pi_\alpha': R\ra R / \gp_\alpha'$, $a\mapsto \gp_\alpha'$. If
$\alpha$ is a limit ordinal then 
\begin{equation}\label{rpaup2}
\gp_\alpha' = \bigcup_{\beta <\alpha } \gp_\alpha'.
\end{equation}

For each  Ore set $S\in \Ore (R)$, let $$\CJ (R, S):=\{ \ga\, | \,
\ga\; {\rm is \; an \; ideal\; of\;} R\; {\rm  such \; that } \;
\pi_\ga (S)\in \Den (R/ \ga , 0)\}$$ where $\pi_\ga : R\ra R/
\ga$, $a\mapsto a+\ga$. It follows from (\ref{paup}) that the
ideals $\gp (S)_l=\bigcup_{\alpha \in \mW}\gp_\alpha$ and $\gp
(S)_r=\bigcup_{\alpha \in \mW}\gp_\alpha'$ coincides since, for
all $\alpha \in \mW$, 
\begin{equation}\label{paup3}
\gp_\alpha = \gp_\alpha'.
\end{equation}
Moreover, for all $\alpha \in \mW$, 
\begin{equation}\label{paup4}
\gp_{\alpha +1}:= \pi_{\alpha}^{-1}  (\ass_l (\pi_\alpha (S),
R/\gp_\alpha ) +\ass_r (\pi_\alpha (S), R/ \gp_\alpha ))\;\; {\rm
and}\;\; \gp (S) = \bigcup_{\alpha \in \mW}\gp_\alpha .
\end{equation}

\begin{theorem}\label{B15Nov10}
Let $S\in \Ore (R)$. Then
\begin{enumerate}
\item $\CJ (R,S) \neq \emptyset$. \item
  $\gp (S)$ is
 the least (with respect to inclusion) element of $\CJ (R, S)$.
 \item If, in addition, $S\in \Den (R, \ga )$ then $\gp (S) =
 \ga$.
\end{enumerate}
\end{theorem}

{\it Proof}. The theorem follows at once from Theorem
\ref{A15Nov10} and its analogue for the right Ore sets in $R$
provided we show that $\gp (S) \neq R$. Suppose that $\gp (S) =R$,
we seek a contradiction. Then $\gp_\alpha \cap S\neq \emptyset$
for some $\alpha$. We can assume that $\alpha$ is the least
possible. Then $\alpha$ is necessarily not a limit ordinal. Then
$0\not\in \pi_{\alpha -1} (S)$ since otherwise we would have
$\gp_{\alpha -1} \cap S\neq \emptyset$, a contradiction.
Therefore, $\pi_{\alpha -1} (S) \in \Ore (R/ \gp_{\alpha -1})$.
Replacing $R$ and $S$ by $R/ \gp_{\alpha -1}$ and $\pi_{\alpha
-1}(S)$ respectively we may assume that $\alpha =1$, i.e. $s'\in
\gp_1=\ass_l (S, R) +\ass_r(S, R)$ for some $s'\in S$. Then $a+b =
s'$ for some elements $a,b\in R$ such that $sa=0$ and $ bt=0$ for
some elements $s,t\in S$. Then $S\ni ss't = s(a+b) t= 0$, a
contradiction. The proof of the theorem is complete. $\Box $


\begin{corollary}\label{T28Nov10}
Let $S$ be an  Ore set in a ring $R$. Then there exists an ordered
pair $(Q, f)$ where $Q$ is a ring and $f: R\ra Q$ is a ring
homomorphism such that

(i)  for all $s\in S$, $f(s)$ is a unit in $Q$; 


 and if $(Q', f')$ is another
 pair satisfying the condition
 (i) then there is a unique ring homomorphism $h: Q\ra Q'$ such that
$f' = hf$. The ring $Q$ is unique up to isomorphism. The ring $Q$
is isomorphic to the  localization of the ring $R/ \gp (S)$ at the
 denominator set $\pi (S)\in \Den (R/ \gp (S), 0)$ where the ideal
$\gp (S)$ of the ring $R$ is defined in (\ref{paup4}) and $\pi :
R\ra R/ \gp (S)$, $ a\mapsto a+\gp (S)$.
\end{corollary}


\section{Examples}\label{EXEM}

In this section, the largest  (left; right; left and right)
 quotient ring  and the maximal (left; right; left and right)
quotient rings are found for the following rings:  the
endomorphism algebra $\End_K(V)$ of an infinite dimensional vector
space with countable basis,  semi-prime Goldie rings,  the algebra
$\mI_1$ of polynomial integro-differential operators, and
Noetherian commutative rings.

{\bf The endomorphism algebra $\End_K(V)$ of an infinite
dimensional vector space $V$ with countable basis}. For a vector
space $V$, let $$\CF = \CF (V) :=\{ \v \in \End_K(V)\, | \,
\dim_K(\ker (\v ))<\infty , \; \dim_K(\coker (\v ))<\infty \}$$ be
the set of {\em Fredholm} linear maps/operators in $V$.

\begin{lemma}\label{a24Nov10}
Let $V$ be an infinite dimensional  vector space with countable
basis and $R:=\End_K(V)$. Then
\begin{enumerate}
\item Let $\v \in R$ and $\dim_K(\im (\v )) =\infty$. Then
\begin{enumerate}
\item $\alpha \v \beta =1$ for some elements $\alpha$ and $\beta$
of $R$ (necessarily, $\alpha$ is a surjection and $ \beta$ is an
injection). Moreover, $\alpha$ and $\beta$ can be chosen to
satisfy the following conditions: $V= \ker (\alpha ) \bigoplus \im
(\v )$ and $V = \ker (\v ) \bigoplus \im (\beta )$. \item If $\v $
is a surjection then $\v  \beta =1$ for some (necessarily,
injective) map $\beta\in R$ that can be chosen to satisfy the
equality $V = \ker (\v ) \bigoplus \im (\beta )$. \item If $\v$ is
an injection then $ \alpha \v =1$ for some (necessarily,
surjective) map $\alpha \in R$ which can be chosen to satisfy the
equality $V = \ker (\alpha ) \bigoplus \im (\v )$.
\end{enumerate}
\item The ideal of compact operators $\CC := \{ \v \in R\, | \,
\dim_K(\im (\v ))<\infty \}$ is the only proper ideal of the ring
$R$. \item Let $ \v \in \CF (V)$ and $c\in R$. Then
\begin{enumerate}
\item If $\v c=0$ then $c\in \CC$. \item If $c\v =0$  then $c\in
\CC$.
\end{enumerate}
\item Let $ \v \in R$. Then
\begin{enumerate}
\item $\ker (\v ) \neq 0$ iff $\v c=0$ for some element $0\neq
c\in \CC$. \item $\im (\v ) \neq V$ iff $c\v =0$ for some element
$ 0\neq c\in \CC$.
\end{enumerate}
\item Let $\v \in R$. Then $\v$ is a right regular element in $R$
(i.e. $\v \psi =0$ implies $ \psi =0$) iff the map $\v : V \ra V$
is an injection. \item Let $\v \in R$. Then $\v$ is a left regular
element in $R$ (i.e. $ \psi\v  =0$ implies $ \psi =0$) iff the map
$\v : V \ra V$ is a surjection. \item Let $\v \in R$. Then $\v$ is
a regular element of $R$ iff $\v$ is a unit of $R$. So,
$\CC_R=\Aut_K(V)$ and $ Q_{cl} (R) = Q_{cl, l} (R) = Q_{cl , r}(R)
= R$.\end{enumerate}
\end{lemma}

{\it Proof}. 1. This is obvious.

2. Statement 2 follows from statement 1.

3a. $\im (c) \subseteq \ker (\v )$, and so $c\in \CC$.

3b. $\ker (c) \supseteq \im (\v )$, and so $\dim_K(\im (\v ))\leq
\dim_K(\coker (\v ))<\infty$, i.e. $c\in \CC$.

4. This is obvious.

5. Statement 5 follows from statement 4a.

6. Statement 6 follows from statement 4b.

7. Statement 7 follows from statements 5 and 6.  $\Box $


\begin{theorem}\label{4Dec10}
Let $V$ be an infinite dimensional vector space with countable
basis and $R:=\End_K(V)$. Then
\begin{enumerate}
\item $\Ass_l(R)= \Ass_r(R) = \Ass (R) = \{ 0, \CC \}$. \item $
S_{l,0}(R) = S_{r,0}(R) = S_0(R) = \Aut_K(V)$ and $Q_l(R) = Q_r(R)
= Q(R) = R$. \item $ S_{l,\CC }(R) = S_{r,\CC }(R) = S_\CC (R) =
\CF $ and $Q_{l, \CC }(R) = Q_{r, \CC }(R) = Q_\CC (R) = R/ \CC $.
\item $\maxAss_l(R) = \maxAss_r(R)= \maxAss (R) = \{ \CC \}$.
\item $R/ \CC $ is a localization maximal ring and a left (resp.
right; left and right) localization  maximal ring.
\end{enumerate}
\end{theorem}

{\it Proof}. 1. Statement 1 follows from statements 2 and 3.

2.  Statement 2 follows from Lemma \ref{a24Nov10}.(7).

4. Statement 4 follows from statement 3.

3. Let $\pi : R\ra R/ \CC$, $ a\mapsto a+\CC$. The ordered pair $(
R/ \CC , \pi )$ satisfies the following conditions

(i) for all $s\in \CF $, $\pi (s)$ is a unit;

(ii) for all $q\in R/ \CC$, $q= \pi ( s)^{-1} \pi ( r)= \pi ( r_1)
\pi (s_1)^{-1}$ for some $s, s_1\in \CF$ and $r,r_1\in R$;

(iii) $\ker ( \pi ) = \CC$ and $\CC = \ass_l(\CF ) = \ass_r (\CF
)$, by Lemma \ref{a24Nov10}.(3).

By Ore's Theorem, $ R/ \CC = \CF^{-1} R= R\CF^{-1}$. To finish the
proof of statement 3 it suffices to show that every regular
element $\pi (r)$ in $R/ \CC$ is invertible.  Since $\pi (r)$
regular in $R/ \CC$, $\dim_K( \ker (r)) <\infty$ (suppose that
$\dim_K( \ker (r ))=\infty$, we seek a contradiction. Fix a
complement subspace $U$ of $\ker ( r )$ in $V$, i.e. $V = \ker ( r
) \bigoplus U$. Let $ p$ be the projection onto $ \ker (r )$. Then
$ \pi ( r) \pi (p)=0$ but $\pi ( p) \neq 0$ since $\dim_K( \ker
(r)) =\infty $, a contradiction).

Similarly, since $\pi ( r)$ is regular in $R/ \CC$, $\dim_K(\coker
( r))<\infty$ (suppose that $\dim_K( \coker (r )) = \infty$, we
seek a contradiction. Fix a complement subspace $W$ of $\im (r )$
in $V$, i.e. $V = \im (r ) \bigoplus W$. Let $q$ be the projection
onto $W$. Then $\pi ( q) \pi ( r) =0$ but $\pi ( q) \neq 0$ since
$ \dim_K( W) = \dim_K(\coker ( r )) = \infty$, a contradiction).
Therefore, $r\in \CF$, and so $\pi (r)$ is a unit.

5. Statement 5 follows from statements 3  and 4, see Theorem
\ref{21Nov10}. $\Box $

$\noindent $

{\bf Semi-prime Goldie rings}.  A ring $R$ is called a {\em left
Goldie} ring if $R$ has finite left uniform dimension  and $R$
satisfies the a.c.c. on left annihilators. Similarly, a right
Goldie ring is defined. A left and right Goldie ring is called a
{\em Goldie} ring. The reader is referred to the books
\cite{Jategaonkar-book}, \cite{MR} and \cite{Stenstrom-RingsQuot}
for more details.

\begin{lemma}\label{a15Dec10}
Let $R$ be a ring.
\begin{enumerate}
\item If $S$ and $T$ are left (resp. right; left and right) Ore
sets in $R$ then so is the semigroup $ST$ generated by $S$ and $T$
in $R$ provided  $0\not\in ST$.  \item If $S\in \Den (R, \ga )$
and $C\in \Den(R,0)$ then $CS\in \Den (R, \ga )$, i.e. the monoid
$CS$ generated by $C$ and $S$ in $R$ is a denominator set with
$\ass (CS ) = \ga$.
\end{enumerate}
\end{lemma}

{\it Proof}. 1. It suffices to prove statement 1, say, for left
Ore sets $S$ and $T$. We have to show that for elements
$s_1t_1s_2t_2\cdots s_nt_n\in ST$ (where  $s_i\in S$ and $t_i\in
T$) and $r\in R$ there exist elements $\th \in ST$ and $r'\in R$
such that $\th r = r's_1t_1 \cdots s_nt_n$. Since $T$ is a left
Ore set, there are elements $t_n'\in T$ and $r_n'\in R$ such that
$t_n'r= r_n't_n$. Since $S$ is a left Ore set, there are elements
$s_n'\in S$ and $r_n''\in R$ such that $s_n'r_n'=r_n''s_n$.
Therefore, $s_n't_n'r= r_n''s_n t_n$. Repeating these two steps
$n-1$ more times we find elements $s_1',\ldots , s_{n-1}'\in S$;
$t_1', \ldots , t_{n-1}'\in T$ and $r_1''\in R$ such that
$s_1't_1'\cdots s_n't_n'r=r_1''s_1t_1\cdots s_nt_n$. Now, set $\th
:= s_1't_1'\cdots s_n't_n'\in S$ and $r':= r_1''\in R$.

2. Let us  show that $CS$ is a multiplicative set. Suppose that
$c_1s_1\cdots c_ns_n=0$ for some elements $c_i\in C$ and $s_i\in
S$, we seek a contradiction. Then $s_1c_2s_2\cdots c_ns_n=0$ since
$c_1$ is a regular element, and so $c_2s_2\cdots c_ns_ns_n'=0$ for
some element $s_n'\in S$ since $S\in \Den (R, \ga )$. Repeating
this argument  several times we come to the situation where
$c's'=0$ for some elements $c'\in C$ and $s'\in S$, i.e. $s'=0$
(since $c'$ is a regular element), a contradiction.

By statement 1,  $CS$ is an Ore set in $R$. To finish the proof of
statement 2  it suffices to show that $\ass (CS)=\ga$, i.e.
$c_1s_1\cdots c_ns_nr =0$ for some $c_i\in C$, $s_i\in S$ and
$r\in R$, implies that $r\in \ga$. The element $c_1$ is a regular
element, hence $s_1\cdots c_ns_nr =0$, and so $c_2s_2\cdots
c_ns_nr s_1'=0$ for some $s_1'\in S$ since $S\in \Den (R , \ga )$.
Repeating the same two steps $n-1$ more times  gives
$rs_1's_2'\cdots s_n'=0$ for some elements $s_i'\in S$, i.e. $r\in
\ga$. $\Box $


\begin{corollary}\label{b15Dec10}
Let $R$ be a prime Goldie ring and $\CC_R$ be the set of regular
elements of $R$. Then $\Ass (R) = \{ 0\}$, $S_0(R) = \CC_R$,
$Q_0(R) = Q_{cl}(R)$ is the only maximal localization of the ring
$R$.
\end{corollary}

{\it Proof}. Let  $\ga \in \Ass (R)$ and $S\in \Den (R, \ga )$. By
Lemma \ref{a15Dec10}.(2), $\CC_R S\subseteq S_\ga (R)$. Since
$Q_{cl}(R)$ is a simple artinian ring (i.e. the matrix ring over a
division ring), we see that $\ga =0$, and so $\CC_R = S_0(R)$.
$\Box $

$\noindent $

Let $R$ be a semi-prime Goldie ring which is not a prime ring and
$\CC_R$ be its set of regular elements. Let $\Min (R) = \{ \gp_1,
\ldots , \gp_s\}$ be the set of minimal primes of the ring $R$. By
Goldie's Theorem, $Q_{cl}(R) := \CC_R^{-1}R \simeq R\CC_R^{-1}
\simeq \prod_{i=1}^s R_i$ is the direct product of simple artinian
rings $R_i$ (i.e. $R_i$ is a matrix ring over a division ring).
The ring $R$ can be identified with its image under the ring
monomorphism  $\s : R\ra Q_{cl}(R)$, $r\mapsto r/1$. For each
$i=1, \ldots , s$, $\gp_i = \s^{-1} (\prod_{i\neq j =1}^s R_j)$
(Proposition 3.2.2, \cite{MR}). For each non-empty set $I$ of the
set $\{ 1, \ldots , s\}$, let the ring homomorphism $\s_I : R \ra
R_I:= \prod_{i\in I}R_i$ be the composition of $\s$ and  the
projection $\prod_{i=1}^sR_i\ra R_I$. Let $R_i^*$ be the group of
units of the ring $R_i$. Then $R_I^*=\prod_{i\in I} R_i^*$ is the
group of units of the ring $R_I$.  A subset $A$ of a set $B$ is
called a {\em proper subset} of $B$ if $A\neq \emptyset, B$.

\begin{theorem}\label{15Dec10}
Let $R$ be a semi-prime Goldie ring which is not a prime ring and
$\{ \gp_1, \ldots , \gp_s\}$ be the set of its minimal prime
ideals.  Then
\begin{enumerate}
\item $\Ass (R) := \{ \ga (I):= \bigcap_{i\in I} \gp_i \, | \,
\emptyset \subsetneqq   I \subseteq \{ 1, \ldots , s\} \}$. \item
For each proper subset $I$ of $ \{ 1, \ldots , s\} $, $S_{\ga
(I)}(R)= \s_{I}^{-1} (R_{I}^*\times R_{CI})$ and $ Q_{\ga (I)} (R)
= Q_{cl}(R)/ R_{CI}\simeq R_{I}$ where $CI:= \{ 1, \ldots ,
s\}\backslash I$. \item $S_0(R) = \CC_R \subseteq S_\ga (R)$ for
all $\ga \in \Ass (R)$. \item $S_{\ga (I)}(R) \subseteq S_{\ga
(J)}(R)$ iff $I\supseteq J$ where $I$ and $J$ are proper subsets
of $\{ 1, \ldots , s\}$. \item If $I$ and $J$ are proper subsets
of $\{ 1, \ldots , s\}$ such that $I\supseteq J$ then, by
statement 4 and the universal property of localization, there is
the unique ring homomorphism $Q_{\ga (I)} (R) = \prod_{i\in
I}R_I\ra Q_{\ga (J)} (R) = \prod_{j\in J}R_j$ which is necessarily
the projection onto $\prod_{j\in J}R_j$ in $\prod_{i\in I}R_i$.
\item $\maxAss (R) = \{ \gp_i, | \, i=1, \ldots , s\}$, $\{
Q_{\gp_i}(R)=R_i \, | \, i=1, \ldots , s\}$ is the set of maximal
localizations of the ring $R$ and $\ass (Q_{\gp_i }) (R) = \gp_i$
for all $i=1, \ldots , s$.
\end{enumerate}
\end{theorem}

{\it Proof}. 1--3. By Theorem \ref{T5Jul10}, $S_0 (R)= \CC_R$ and
$Q(R) = Q_{cl}(R)$. By Lemma \ref{a15Dec10}.(2), $\CC_R \subseteq
S_\ga (R)$ for all $\ga \in \Ass (R)$ (and statement 3 follows),
and there is a natural homomorphism $Q_{cl } (R) \ra  Q_\ga (R)$,
$r\mapsto r/1$.  By Lemma \ref{a5Jul10}, the ring $Q_\ga (R)$ is a
localization of the ring $Q_{cl}(R)$. Moreover, by Proposition
\ref{TA21Nov10}.(2), the map $\Ass (R) \ra \Ass (Q_{cl}(R))$, $\ga
\mapsto S^{-1}\ga$, is a bijection with the inverse map
\begin{equation}\label{TCC1}
\gb \mapsto R\cap \gb ,
\end{equation}
and the map $$ \{ S_\ga (R)\, | \, \ga \in \Ass  (R) \} \ra \{
S_\gb (Q_{cl}(R))\, | \, \gb \in \Ass (Q_{cl}(R))\}, \;\; S_{R\cap
\gb } (R) \mapsto S_\gb (Q_{cl}(R)), $$ is a bijection and
\begin{equation}\label{TCC2}
S_{R\cap \gb } (R) = \s^{-1}_{R\cap \gb } (S_\gb (Q_{cl }(R)))
\end{equation}
where 
\begin{equation}\label{TCC3}
\s_{R\cap \gb } : R\ra Q_{R\cap \gb } (R) = S_{R\cap \gb }
(R)^{-1}R \simeq S_\gb (Q_{cl}(R))^{-1}Q_{cl}(R).
\end{equation}
Clearly, $\Ass (Q_{cl}(R)) = \Ass (\prod_{i=1}^s R_i) = \{ 0\}
\bigcup \{ R_I, | \, \emptyset  \neq I \varsubsetneqq \{ 1, \ldots
, s\} \}$ and $R_I = Q_{cl}(R)/ R_{CI}= Q_{R_{CI}}(Q_{cl}(R))$ (by
Theorem \ref{T25Nov10}, since every regular element of $R_I$ is a
unit). Now, statement 1 follows from (\ref{TCC1}) and the fact
that $\gp_i = R\cap \prod_{j\neq i} R_j$. For every proper subset
$I$ of $\{ 1, \ldots ,s\}$, $R\cap R_I = \ga (CI)$, $\s_{R\cap
R_I} = \s_{CI}$ and $S_{R_I}(Q_{cl}(R)) = R^*_{CI}\times R_I$.
Then, by (\ref{TCC2}), $S_{R\cap R_I}(R) = \s^{-1}_{R\cap R_I}
(S_{R_I}(Q_{cl}(R)))= \s^{-1}_{CI}(R_{CI}^*)$, and, by
(\ref{TCC3}), $Q_{\ga (CI)}(R)=Q_{R\cap R_I}(R) \simeq
S_{R_I}(Q_{cl}(R))^{-1} Q_{cl}(R) \simeq R_{CI}$.

4. Statement 4 follows from statement 2.

5. Statement 5 follows from statement 2.

6. Statement 6 follows from statements 2, 4 and 5.  $\Box $

$\noindent $

{\bf The algebra $\mI_1$ of polynomial integro-differential
operators}. Let us collect some facts for the algebra $\mI_1$
which are used in the proofs of Theorem \ref{A4Dec10} and
Proposition \ref{a11Dec10}. For the details the reader is referred
to \cite{algintdif} or \cite{intdifline}. Throughout,  $K$ is a
field of characteristic zero and  $K^*$ is its group of units;
$P_1:= K[x]$ is a polynomial algebra in one variable $x$ over $K$;
$\der:=\frac{d}{d x}$; $\End_K(P_1)$ is the algebra of all
$K$-linear maps from $P_1$ to $P_1$,  and $\Aut_K(P_1)$ is its
group of units (i.e. the
 group of all the invertible linear maps from $P_1$ to $P_1$); the
subalgebras  $A_1:= K \langle x , \der \rangle$ and
 $\mI_1:=K\langle x,\der ,  \int\rangle $ of $\End_K(P_1)$
  are called the (first) {\em Weyl algebra} and the {\em algebra  of polynomial
integro-differential operators} respectively  where $\int: P_1\ra
P_1$, $ p\mapsto \int p \, dx$, is the  {\em integration},  i.e.
$\int : x^n \mapsto \frac{x^{n+1}}{n+1}$ for all $n\in \N$.  The
algebra $\mI_1$ is neither left nor right Noetherian and not a
domain. Moreover, it contains infinite direct sums of nonzero left
and right ideals, \cite{algintdif}.

The algebra $\mI_1$  is generated by the elements $\der $, $H:=
\der x$ and $\int$ (since $x=\int H$) that satisfy the defining
relations (Proposition 2.2, \cite{algintdif}): $$\der \int = 1,
\;\; [H, \int ] = \int, \;\; [H, \der ] =-\der , \;\; H(1-\int\der
) =(1-\int\der ) H = 1-\int\der .$$ The elements of the algebra
$\mI_1$,  
\begin{equation}\label{eijdef}
e_{ij}:=\int^i\der^j-\int^{i+1}\der^{j+1}, \;\; i,j\in \N ,
\end{equation}
satisfy the relations $e_{ij}e_{kl}=\d_{jk}e_{il}$ where $\d_{jk}$
is the Kronecker delta function. Notice that
$e_{ij}=\int^ie_{00}\der^j$.

The algebra $\mI_1=\bigoplus_{i\in \Z} \mI_{1, i}$ is a
$\Z$-graded algebra ($\mI_{1, i} \mI_{1, j}\subseteq \mI_{1, i+j}$
for all $i,j\in
\Z$) where 
\begin{equation}\label{I1iZ}
\mI_{1, i} =\begin{cases}
D_1\int^i=\int^iD_1& \text{if } i>0,\\
D_1& \text{if }i=0,\\
\der^{|i|}D_1=D_1\der^{|i|}& \text{if }i<0,\\
\end{cases}
\end{equation}
 the algebra $D_1:= K[H]\bigoplus \bigoplus_{i\in \N} Ke_{ii}$ is
a commutative non-Noetherian subalgebra of $\mI_1$, $ He_{ii} =
e_{ii}H= (i+1)e_{ii}$  for $i\in \N $ (notice that
$\bigoplus_{i\in \N} Ke_{ii}$ is the direct  sum of non-zero
ideals of $D_1$); $(\int^iD_1)_{D_1}\simeq D_1$, $\int^id\mapsto
d$; ${}_{D_1}(D_1\der^i) \simeq D_1$, $d\der^i\mapsto d$,   for
all $i\geq 0$ since $\der^i\int^i=1$.
 Notice that the maps $\cdot\int^i : D_1\ra D_1\int^i$, $d\mapsto
d\int^i$,  and $\der^i \cdot : D_1\ra \der^iD_1$, $d\mapsto
\der^id$, have the same kernel $\bigoplus_{j=0}^{i-1}Ke_{jj}$.

Each element $a$ of the algebra $\mI_1$ is the unique finite sum
\begin{equation}\label{acan}
a=\sum_{i>0} a_{-i}\der^i+a_0+\sum_{i>0}\int^ia_i +\sum_{i,j\in
\N} \l_{ij} e_{ij}
\end{equation}
where $a_k\in K[H]$ and $\l_{ij}\in K$. This is the {\em canonical
form} of the polynomial integro-differential operator
\cite{algintdif}. The algebra $\mI_1$ has the only proper ideal
$F=\bigoplus_{i,j\in \N}Ke_{ij} \simeq M_\infty (K)$ and
 $F^2= F$. The factor algebra $\mI_1/F$ is canonically isomorphic to
the skew Laurent polynomial algebra $B_1:= K[H][\der, \der^{-1} ;
\tau ]$, $\tau (H) = H+1$, via $\der \mapsto \der$, $ \int\mapsto
\der^{-1}$, $H\mapsto H$ (where $\der^{\pm 1}\alpha = \tau^{\pm
1}(\alpha ) \der^{\pm 1}$ for all elements $\alpha \in K[H]$). The
algebra $B_1$ is canonically isomorphic to the (left and right)
localization $A_{1, \der }$ of the Weyl algebra $A_1$ at the
powers of the element $\der$ (notice that $x= \der^{-1} H$).
Therefore, they have the common skew field of fractions, $\Frac
(A_1) = \Frac (B_1)$, the {\em first Weyl skew field}.

 The algebra
$\mI_1$ admits the involution $*$ over the field $K$: $\der^* =
\int$, $\int^* = \der$ and $H^* = H$, i.e. it is a $K$-algebra
{\em anti-isomorphism} ($(ab)^* = b^* a^*$) such that $a^{**} =a$.
Therefore, the algebra $\mI_1$ is {\em self-dual}, i.e. it is
isomorphic to its opposite algebra $\mI_1^{op}$. As a result, the
left and right properties of the algebra $\mI_1$ are the same.
Clearly, $e_{ij}^* = e_{ji}$ for all $i,j\in \N$, and so $F^* =
F$.

The next theorem describes the one-sided largest quotient rings of
the algebra $\mI_1$.
\begin{theorem}\label{CC18Jun10}
{\rm (Theorem 9.7, \cite{intdifline})}
\begin{enumerate}
\item $S_{r,0}(\mI_1) = \mI_1 \bigcap \Aut_K(K[x])$ and the
largest regular right quotient ring $Q_r(\mI_1)$ of $\mI_1$ is the
subalgebra of $\End_K(K[x])$ generated by $\mI_1$ and $S_{r,0}
(\mI_1)^{-1}:= \{ s^{-1}\, | \, s\in  S_{r,0} (\mI_1)\}$.  \item
$S_{l,0}(\mI_1) =S_{r,0}(\mI_1)^*$ and $S_{l,0}(\mI_1) \neq
S_{r,0}(\mI_1)$. \item The rings $Q_l(\mI_1)$ and $Q_r(\mI_1)$ are
not isomorphic.
\end{enumerate}
\end{theorem}

The next theorem describes the (two-sided) largest quotient ring
of the algebra $\mI_1$, it is tiny comparing with the one-sided
largest quotient rings of the algebra $\mI_1$.
\begin{theorem}\label{A4Dec10}
Let $\CM := (K[H]+F) \bigcap \Aut_K( K [x])$. Then
\begin{enumerate}
\item $S_0( \mI_1) = S_{l,0}(\mI_1) \bigcap S_{r,0}(\mI_1) $,
$S_0( \mI_1)$ is a proper subset of the sets  $S_{l,0}(\mI_1)$ and
$S_{r,0}(\mI_1) $, and $S_0(\mI_1)^* = S_0(\mI_1)$ where $*$ is
the involution of the algebra $\mI_1$. \item $ S_{l,0}(\mI_1)
\bigcap S_{r,0}(\mI_1) =\CM $ and $\CM$ is the set of regular
elements of the algebra $K[H]+ F$. \item Let $\CM_0:= D_1\bigcap
\Aut_K(K[x])$. Then $\CM_0\subseteq \CM$, $\CM = \CM_0(1+F)^*=
(1+F)^* \CM_0$ and $ \CM_0\bigcap (1+F)^* = (1+F_0)^*$ where
$F_0:= \bigoplus_{i\in \N}Ke_{ii}$. \item $\CM_0$ is the set of
regular elements of the commutative (non-Noetherian)
 algebra $D_1$; $D_1= \CM_0 (1+F_0) \coprod F_0= \CM_0\cup \{ 0\} +F_0$, $Q_{cl}
 (D_1) := \CM_0^{-1}D_1 = \CM_0^{-1}\CM_0(1+F_0)\coprod F_0 =
 \CM_0^{-1}\CM_0\cup \{ 0\} +F_0$.
 \item $Q(\mI_1) = S_0(\mI_1)^{-1} \mI_1 = \sum_{i\in \Z}
 Q_{cl}(D_1)v_i +F= \sum_{i\in \Z } (\CM_0^{-1}\CM_0 \cup \{ 0\} )v_i +F=
 \sum_{i\in \Z} v_i Q_{cl } (D_1) +F = \sum_{i\in \Z} v_i
 (\CM_0^{-1}\CM_0\cup \{ 0\} ) +F$ where $Q_{cl } (D_1)$ is the classical ring of
 fractions of the commutative ring $D_1$ and
 $$ v_i:= \begin{cases}
\int^i& \text{if }i\geq 1,\\
1& \text{if }i=0,\\
\der^{|i|}&\text{if }, i\leq -1.
\end{cases}$$
\item $Q(\mI_1) \subsetneqq Q_l(\mI_1)$ and $Q(\mI_1) \subsetneqq
 Q_r(\mI_1)$.
\end{enumerate}
\end{theorem}
 {\it Proof}. 2. Recall that $S_{r,0} (\mI_1) = \mI_1 \bigcap
 \Aut_K( K[x])$ (Theorem 9.7.(2), \cite{intdifline}),
 $S_{l,0} (\mI_1) =S_{r,0}(\mI_1)^*$ (Theorem 9.7.(3),
 \cite{intdifline}), and if $a\in \mI_1\backslash F$ and $ a\in
 \mI_1\bigcap \Aut_K( K[x])$ then $ a\in \sum_{i\leq 0} D_1v_i+F$
 (Theorem 6.2.(1), \cite{intdifline}). Since $ (v_i)^* = v_{-i}$
 and $d^* = d$ for all elements $ d\in D_1$,
\begin{eqnarray*}
S_{l,0}(\mI_1) \bigcap
 S_{r,0}( \mI_1) &=&S_{l,0}(\mI_1) \bigcap
 S_{l,0}( \mI_1)^* \subseteq (K[H]+F)\bigcap \Aut_K(K[x])
  \\
 &\subseteq & S_{l,0}(\mI_1)\bigcap
 S_{l,0}( \mI_1)^* =  S_{l,0}(\mI_1) \bigcap
 S_{r,0}( \mI_1)\\
\end{eqnarray*}
 i.e. $S_{l,0}(\mI_1) \bigcap  S_{r,0}( \mI_1)=\CM$.

 Let $\CR $ be the set of regular elements of the algebra
 $K[H]+F$. Clearly, $\CM \subseteq \CR$. Let $r\in \CR$. Then
 $r\not\in F$. By Proposition 6.1.(1), \cite{intdifline},
 $\ind_{K[x]}(r) =0$. Then, by Theorem 6.2.(3), \cite{intdifline},
 $r_{K[x]}$ is a bijection iff $r_{K[x]}$ is an injection iff
 $r_{K[x]}$ is a surjection where $r_{K[x]}: K[x]\ra K[x]$, $ p\mapsto rp$. Therefore, to prove the equality $\CM
 =\CR$ it suffices to show that $r_{K[x]}$ is an injection.
 Suppose that $\ker (r_{K[x]})\neq 0$, we seek a contradiction.
 Since $r\not\in F$, by Lemma 6.10, \cite{intdifline}, there
 exists an idempotent $f\in F$ such that $\ker (r_{K[x]}) = \im
 (f_{K[x]})$. In particular, $rf=0$ but $f\neq 0$, a
 contradiction. Therefore, $\CM = \CR$.

 1. Since $S_{l,0} (\mI_1) \neq S_{r,0}(\mI_1)$ (Theorem
 9.7.(3,4), \cite{intdifline}), the set $S_0(\mI_1)$ is a proper
 subset of the sets $S_{l,0} (\mI_1)$ and $ S_{r,0}(\mI_1)$.
 Clearly,
 \begin{eqnarray*}
 S_0(\mI_1) &\subseteq & \CM =S_{l,0} (\mI_1) \bigcap  S_{r,0}(\mI_1)=
 S_{l,0} (\mI_1) \bigcap  S_{r,0}(\mI_1)^* \\
 &=&(S_{l,0} (\mI_1) \bigcap  S_{l,0}(\mI_1)^*)^*=
 S_{l,0} (\mI_1) \bigcap  S_{r,0}(\mI_1)=\CM .
\end{eqnarray*}
We have used the fact that $ S_{l,0} (\mI_1) = S_{r,0}(\mI_1)^*$
(Theorem 9.7.(3), \cite{intdifline}). In view of statement 2 and
the equality $\CM^* = \CM$, to finish the proof of statement 1 it
suffices to show that $\CM $ is a left Ore set in $\mI_1$,  that
is, for any elements $a\in \mI_1$ and $u= \alpha +u_1\in \CM$
where $\alpha \in K[H]$ and $u_1\in F$, there exists an element
$b\in \CM$ such that $bau^{-1} \in \mI_1$ (where the product
$bau^{-1}$ is taken in $\End_K(K[x])$, recall that $\mI_1\subseteq
\End_K(K[x])$). The element $a$ is a sum $
\sum_{i=-n}^n\alpha_iv_i +f$ for some natural number $n$ where $
\alpha_i\in K[H]$.

For all elements $\beta \in K[H]$, $\beta v_i = v_i \tau^i (\beta
)$ where $\tau \in \Aut_K( K[H])$ and  $\tau (H) = H+1$. Notice
that $F\CM^{-1} \subseteq F$ and $ \alpha \cdot u^{-1} = \alpha (
\alpha +u_1)^{-1} \in (1+F)^*\subseteq \CM$. Let $b':=
\prod_{i=-n}^n \tau^i (\alpha )$. Then
$$ b'au^{-1}=  b'(\sum_{i=-n}^n \alpha_i v_iu^{-1} +fu^{-1})
= \sum_{i=-n}^n \alpha_iv_i\frac{\tau^i(b')}{\alpha}\alpha u^{-1}
+b'fu^{-1}\in \mI_1,$$ since $\frac{\tau^i(b')}{\alpha}\in K[H]$,
$\alpha u^{-1} \in \CM$ and  $b'fu^{-1} \in b'\cdot F \CM^{-1}
\subseteq b' \cdot F \subseteq F$. Let $I:=\{ i\in \N \, | \,
b*x^i=0\}$. Since $H*x^i = (i+1)x^i$ for all $i\in \N$, we see
that $I=\{ i\in \N \, | \, i+1$ is a root of the polynomial $b'\in
K[H]\}$. Then the element $b:= b'+\sum_{i\in I}e_{ii}\in \CM$ and
$ba u^{-1}\in \mI_1$.

3. Since $D_1= K[H]+F_0\subseteq K[H]+F$, the inclusion
$\CM_0\subseteq \CM$ follows. It is obvious that $\CM_0\bigcap
(1+F)^* = (1+F_0)^*$. To prove the equality $ \CM = \CM_0(1+F)^*$
we have to  show that each element $u\in \CM$ is a product $vw$
for some elements $v\in \CM_0$ and $w\in (1+F)^*$. Notice that $u
= \alpha +f$ for some elements $\alpha \in K[H]$ and $f\in F$.
Choose an element, say $g\in F_0$, such that $v:= \alpha +g\in
\CM_0$ (see the proof of statement 1). Then $u= \alpha +g+f-g= v(
1+v^{-1} (f-g)) = vw$ where $w:= 1+v^{-1} (f-g) \in (1+F)^*$ since
$\CM_0^{-1} F\subseteq F$. Therefore, $\CM = \CM_0 (1+F)^*$. Then
$\CM = (1+F)^* \CM_0$ since, for all elements $s\in \CM_0$,
$s(1+F)^*s^{-1} \subseteq (1+F)^*$.

4. Let $\CR$ be the set of regular elements of the ring $D_1$.
Then $ \CM_0\subseteq \CR$. Let $r\in \CR$, in order to prove that
the inverse inclusion holds, $ \CM_0\supseteq \CR$,   we have to
show that $r\in \CM_0$ or, equivalently, that $r\in \Aut_K(
K[x])$. Notice that $K[x]=\bigoplus_{i\in \N} Kx^i$ is the direct
sum of $r$-invariant subspaces, i.e. $r* Kx^i \subseteq Kx^i$ for
all $i\in \N$. If $r\not\in \Aut_K(K[x])$ then $r*Kx^i =0$ for
some $i$, and so $re_{ii} =0$ since $\im ( e_{ii})=Kx^i$ where
$e_{ii}\in D_1$, a contradiction. Therefore, $r\in \Aut_k( K[x])$,
and so $\CM_0 = \CR$.

 Clearly, $D_1\supseteq
\CM_0(1+F_0) \coprod F_0$. To prove that the reverse inclusion
holds we have to show that every element $d\in D_1\backslash F_0$
belongs to the set $\CM_0 (1+F_0)$. Notice that $d = \alpha +f$
for some elements $0\neq \alpha \in K[H]$ and $f\in F_0$. Choose
an element, say $ g\in F_0$, such that $u:= \alpha +g \in \CM_0$.
Then $d= \alpha +g+f-g = u( 1+u^{-1} (f-g)) \in \CM_0(1+F_0)$
since $\CM_0^{-1} F_0\subseteq F_0$. Therefore, $D_1= \CM_0
(1+F_0) \coprod F_0$, and so $D_1= \CM_0\cup \{ 0\} +F_0$. Then
$\CM_0^{-1} D_1 = \CM_0^{-1}\CM (1+F) \coprod F_0=
\CM_0^{-1}\CM\cup \{ 0\} +F_0$ since $\CM^{-1} F_0\subseteq F_0$.

5. By statements 3 and 4, for all $i\in \Z$,
\begin{eqnarray*}
\CM^{-1} D_1v_i &=& \CM_0^{-1} (1+F)^* D_1v_i = \CM_0^{-1} (1+F)^*
(\CM_0\cup \{ 0\} +F_0)v_i \\
&\subseteq &  (\CM_0^{-1} \CM_0\cup \{ 0\} +F)v_i \subseteq
(\CM_0^{-1}\CM_0\cup \{ 0\} )v_i+F.
\end{eqnarray*}
 Therefore,
\begin{eqnarray*}
 S_0(\mI_1)^{-1}\mI_1&=& \CM^{-1}\mI_1=\sum_{i\in \Z}\CM^{-1}D_1v_i+\CM^{-1}F = \sum_{i\in \Z}(\CM_0^{-1}\CM_0\cup \{ 0\} )v_i+F\\
 &=& \sum_{i\in \Z} (\CM_0^{-1}\CM_0\cup \{ 0\} +F_0) v_i+F = \sum_{i\in \Z} Q_{cl} (D_1) v_i +F. \\
\end{eqnarray*}
Applying the involution $*$ to the inclusion $\CM^{-1}D_1v_i
\subseteq (\CM_0^{-1}\CM_0\cup \{ 0\} )v_i +F$  and using the
equalities $\CM^* = \CM$, $D_1^* = D_1$ and $v_i^* = v_{-i}$ we
obtains the inclusion $v_{-i} D_1\CM^{-1} \subseteq
v_{-i}(\CM_0^{-1}\CM_0\cup \{ 0\} ) +F$ for all $i\in \Z$. Then
\begin{eqnarray*}
 \mI_1S_0(\mI_1)^{-1} &=& \mI_1\CM^{-1} = \sum_{i\in \Z} v_iD_1\CM^{-1} +F\CM^{-1} = \sum_{i\in \Z} v_i (\CM_0^{-1}\CM_0\cup \{ 0\} )+F\\
 &=& \sum_{i\in \Z} v_i (\CM_0^{-1}\CM_0\cup \{ 0\} +F)+F= \sum_{i\in \Z} v_i Q_{cl} (D_1) +F. \\
\end{eqnarray*}
6. The inclusion $Q(\mI_1) \subseteq Q_l( \mI_1)$ (resp. $
Q(\mI_1) \subseteq Q_r(\mI_1)$) is a proper inclusion by statement
5 and Theorem 9.7, \cite{intdifline}. Another way to prove this is
as follows. Suppose that $Q(\mI_1) = Q_l(\mI_1)$, we seek a
contradiction. The algebra $\mI_1$ has the involution $*$, and so
it is isomorphic to its opposite algebra of $\mI_1$.  Therefore,
$Q(\mI_1) = Q_r(\mI_1)$, but the rings $Q_l(\mI_1)$ and
$Q_r(\mI_1)$ are not isomorphic, by  Theorem 9.7.(4),
\cite{intdifline}, a contradiction. $\Box $


\begin{proposition}\label{a11Dec10}
\begin{enumerate}
\item $\Ass_l(\mI_1) = \Ass_r (\mI_1) = \Ass (\mI_1)= \{ 0, F\}$
and $\maxAss_l(\mI_1) = \maxAss_r(\mI_1) =\maxAss (\mI_1) = \{
F\}$. \item $S_{l,F} (\mI_1) = S_{r, F} (\mI_1) =S_{F}
(\mI_1)=\mI_1\backslash F$ and $ Q_{l,F}(\mI_1) = Q_{r, F}(\mI_1)
= Q_F(\mI_1) =\Frac (B_1) = \Frac (A_1)$. \item $\maxDen_l(\mI_1)
= \maxDen_r(\mI_1) = \maxDen (\mI_1) =\{ \mI_1\backslash F\}$.
\end{enumerate}
\end{proposition}

{\it Proof}. 2. Let $\pi_F : \mI_1\ra \mI_1/F = B_1$, $a\mapsto
a+F$. Then $\pi_F^{-1} (B_1\backslash \{ 0\} ) = \mI_1\backslash
F$ is a multiplicative set of the algebra $\mI_1$ with
$\ass_l(\mI_1\backslash F) = \ass_r(\mI_1\backslash F) =F$ since
for all $i,j\in \N$, $\der^{i+1}e_{ij}=0$, $ e_{ij}\int^{i+1} =0$
and $F=\bigoplus_{i,j\in \N} Ke_{ij}$ is the only proper ideal of
the algebra $\mI_1$ and the algebra $\mI_1/F = B_1$ is a domain.
Since $B_1$ is a Noetherian domain and $\Frac (A_1) = \Frac (B_1)=
(B_1\backslash \{ 0\} )^{-1} B_1 = B_1(B_1\backslash \{ 0\} )^{-1}
$, we see that, by the universal property of localization, $\Frac
(\mI_1) = (\mI_1\backslash F)^{-1} \mI_1= \mI_1(\mI_1\backslash
F)^{-1} $ and so $\mI_1\backslash F$ is a left and right
denominators set of the algebra $\mI_1$ with $\ass
(\mI_1\backslash F) = F$. Now, statement 2 follows from Theorem
\ref{15Nov10}.(1).

1. Statement 1 follows from statement 2.

3. Statement 3 follows from statements 1 and 2.  $\Box $

$\noindent $

{\bf Noetherian commutative rings}.

{\it Example}. If $R$ is a commutative ring with a single minimal
prime ideal $\gp$ (eg, $R$ is a commutative domain) then $\Ass (R)
= \{ 0\}$, $S_0(R) = R\backslash \gp$ and $Q(R) = R_\gp$, the
localization of the ring $R$ at the prime ideal $\gp$.

\begin{lemma}\label{a12Dec10}
Let $R$ be a commutative ring, $\CC_R$ be its set of regular
elements and $S\in \Ore (R)$. Then
\begin{enumerate}
\item $\gp (S) = \ass (S)$. \item $\CC_R S\in \Den (R, \ass (S))$.
\end{enumerate}
\end{lemma}

{\it Proof}. 1. The equality follows from the definition of the
ideal $\gp (S)$, see (\ref{paup3}) and  (\ref{paup4}).

2. Obvious. $\Box $

$\noindent $

Let $R$ be a Noetherian commutative ring and $\Min (R) = \{ \gp_1,
\ldots , \gp_s\}$ be its set of minimal primes with $s\geq 2$.
Then $S_0= S_0(R) = R\backslash \bigcup_{i=1}^s \gp_i$ is the set
of regular elements of the ring $R$ (i.e. the set of
non-zero-divisors of $R$) and $S_0^{-1}R\simeq \prod_{i=1}^s R_i$
is the direct product of local commutative artinian rings $(R_i,
\gm_i)$ where $\gm_i$  is the maximal ideal of the ring
$R_i=R_{\gp_i}$, the localization of the ring $R$ at $\gp_i$. For
each non-empty subset $I$ of the  set $\{ 1, \ldots, s\}$, let the
ring homomorphism $\s_I : R\ra R_I := \prod_{i\in I} R_i$ be the
composition of the ring monomorphism $\s : R\ra S_0^{-1}R$,
$r\mapsto r/1$ and the projection $S_0^{-1}R\ra R_I$.  The group
of units $R_I^*$ of the ring $R_I$ is equal to $\prod_{i\in I}
R_i^*$ where $R_i^*= R_i\backslash \gm_i$ is the group of units of
the ring $R_i$. For each number $i=1, \ldots , s$,  $\gp_i'
:=\gm_i\times  \prod_{j\neq i} R_j$ is an ideal of the ring
$S_0^{-1}R$. Then $\gp_i'': = \s^{-1} (\gp_i')$ is an ideal of the
ring $R$.

\begin{proposition}\label{12Dec10}
Let $R$ be a Noetherian commutative  ring, $\{ \gp_1, \ldots ,
\gp_s\}$ be the set of its minimal prime ideals, $s\geq 2$ and $\{
\gp_1'', \ldots , \gp_s''\}$ be as above. Then
\begin{enumerate}
\item $\Ass (R) :=\{ 0\} \bigcup  \{ \ga (I):= \bigcap_{i\in I}
\gp_i'' \, | \, \emptyset \subsetneqq   I \subsetneqq  \{ 1,
\ldots , s\} \}$. \item For each proper subset $I$ of $ \{ 1,
\ldots , s\} $, $S_{\ga (I)}(R)= \s_{I}^{-1} (R_{I}^*\times
R_{CI})=R\backslash \bigcup_{i\in I}\gp_i$ and $ Q_{\ga (I)} (R)
\simeq R_{I}$ where $CI:= \{ 1, \ldots , s\}\backslash I$. \item
$S_0(R) = R\backslash \bigcup_{i=1}^s \gp_i \subseteq S_\ga (R)$
for all $\ga \in \Ass (R)$. \item $S_{\ga (I)}(R) \subseteq S_{\ga
(I)}(R)$ iff $I\supseteq J$ where $I$ and $J$ are proper subsets
of $\{ 1, \ldots , s\}$. \item If $I$ and $J$ are proper subsets
of $\{ 1, \ldots , s\}$ such that $I\supseteq J$ then, by
statement 4 and the universal property of localization, there is
the unique ring homomorphism $Q_{\ga (I)} (R) = \prod_{i\in
I}R_I\ra Q_{\ga (J)} (R) = \prod_{j\in J}R_j$ which is necessarily
the projection onto $\prod_{j\in J}R_j$ in $\prod_{i\in I}R_i$.
\item $\maxAss (R) = \{ \gp_i'', | \, i=1, \ldots , s\}$, $\{
Q_{\gp_i''}(R)=R_{\gp_i} \, | \, i=1, \ldots , s\}$ is the set of
maximal localizations of the ring $R$.
\end{enumerate}
\end{proposition}

{\it Proof}. 1--3. Clearly, $S_0:= S_0 (R)= \CC_R$ where $\CC_R
=R\backslash  \bigcup_{i=1}^s \gp_i$ is the set of regular
elements of the ring $R$. By Lemma \ref{a12Dec10}.(2), $S_0
\subseteq S_\ga (R)$ for all $\ga \in \Ass (R)$ and statement 3
follows. We identify the ring $R$ with its image  under the ring
monomorphisms $R\ra S_0^{-1}R = \prod_{i=1}^s R_i$, $r\mapsto
r/1$. By statement 3 and Proposition \ref{A21Nov10}.(2), the map
$\Ass (R) \ra \Ass (S_0^{-1}R)$, $\ga \mapsto S^{-1}\ga$, is a
bijection with the inverse map 
\begin{equation}\label{CC1}
\gb \mapsto R\cap \gb ,
\end{equation}
and the map $$ \{ S_\ga (R)\, | \, \ga \in \Ass  (R) \} \ra \{
S_\gb (S_0^{-1}R)\, | \, \gb \in \Ass (S_0^{-1}R)\}, \;\; S_{R\cap
\gb } (R) \mapsto S_\gb (S_0^{-1}R), $$ is a bijection and
\begin{equation}\label{CC2}
S_{R\cap \gb } (R) = \s^{-1}_{R\cap \gb } (S_\gb (S_0^{-1}R))
\end{equation}
where 
\begin{equation}\label{CC3}
\s_{R\cap \gb } : R\ra Q_{R\cap \gb } (R) = S_{R\cap \gb }
(R)^{-1}R \simeq S_\gb (S_0^{-1}R)^{-1}S_0^{-1}R.
\end{equation}
{\it Claim:} $\Ass (S_0^{-1}R) = \{ 0 \} \bigcup \{ R_I \, | \,
\emptyset \varsubsetneqq I \varsubsetneqq \{ 1, \ldots , s\}\}$.
Let $\CR$ be the RHS of the equality. Then $\Ass (S_0^{-1}R)
\supseteq \CR$ since $R_I = \ass (\s_{CI}^{-1} (R_{CI}^*))$ where
$CI$ is the complement of the proper set $I$ in $\{ 1, \ldots ,
s\}$. To prove that the reverse inclusion holds, i.e. $\Ass
(S_0^{-1}R)\subseteq \CR$, we have to show that, for each
multiplicatively closed subset $S$ of $S^{-1}_0R$ that does not
entirely consists of non-zero-divisors, $\ass (S) = R_I$ for some
proper subset $I$ of $\{ 1, \ldots , s\}$. For each element $r\in
R$, the subset of $\{ 1, \ldots , s\}$, $\supp (r) :=\{ i \, | \,
\pi_i(r) \in R_i^*$ and $ i\in \{ 1, \ldots , s\} \}$ is called
the {\em support} of the element $r$ where $\pi_i : R\ra S_0^{-1}R
= \prod_{j=1}^s R_j\ra R_i$. For all elements $s,t\in S$, $\supp
(st) = \supp (s) \cap \supp (t)$. Let $\supp (S):= \bigcap_{s\in
S} \supp (s)$. Clearly, there exists an element, say $s\in S$,
such that $\supp (s) = \Supp (S)$ (eg, $s$ is the product of all
elements of $S$ with all possible distinct supports).  Replacing
$s$ with $s^n$ for some natural number we may additionally assume
that $\pi_i(s)=0$ for all elements $j\not\in I$ (take $n$ such
that $\gm_i^n =0$ for all $i=1, \ldots , s$). Then clearly $\ass
(S) = \ass (\{ s^i\}_{i\in \N}) = R_{C\supp (S)}$. This finishes
the proof of the Claim.

Statement 1 follows from the Claim and (\ref{CC1}). For every
proper subset $I$ of $\{ 1, \ldots , s\}$, $R\cap R_I = \ga (CI)$,
$\s_I = \s_{R\cap R_{CI}}$ and $S_{R_I}(S_0^{-1}R)=R_{CI}^*\times
R_I$. Then, by (\ref{CC2}), $S_{R\cap R_I}(R) = \s^{-1}_{R\cap
R_I}(S_{R_I}(S_0^{-1}R))= \s^{-1}_{CI}(R^*_{CI}\times R_I)$ and,
by (\ref{CC3}), $Q_{R\cap R_I}(R) \simeq
S_{R_I}(S_0^{-1}R)^{-1}(S_0^{-1}R)\simeq R_{CI}$. To finish the
proof of statement 2 it remains to show that $S_{R\cap R_I} (R) =
R\backslash \bigcup_{i\in I}\gp_i$ for each  proper subset $I$ of
the set $\{ 1, \ldots , s\}$. An element $s\in R$ belongs to the
set $S_{\ga (CI)}=S_{R\cap R_I} (R) =\s^{-1}_{CI}(R^*_{CI}\times
R_I)$ iff $\s (r) \in R^*_{CI}\times R_I$ iff $s\not\in
\bigcup_{i\in CI} \gp_i$ iff $s\in R\backslash \bigcup_{i\in CI}
\gp_i$.

4. Statement 4 follows from statement 2.

5. Statement 5 follows from statement 2.

6. Statement 6 follows from statements 2, 4 and 5.  $\Box $

$\noindent $

Department of Pure Mathematics

University of Sheffield

Hicks Building

Sheffield S3 7RH

UK

email: v.bavula@sheffield.ac.uk

\end{document}